\documentclass[10pt,a4paper]{amsart}
\usepackage[utf8]{inputenc}
\usepackage{amsmath,amsfonts,amssymb,amsthm}
\usepackage[colorlinks=true, pdfstartview=FitV, linkcolor=blue, citecolor=blue, urlcolor=blue,pagebackref=false]{hyperref}
\usepackage{macros}

\usepackage[margin=1.5in]{geometry}
\renewcommand*{\hat}[1]{\widehat{#1}}

\renewcommand*{\bar}[1]{\overline{#1}}

\newcommand{\ind}[1]{\textbf{1}_{\{#1\}}}

\newcommand{\ou}{\<1_black> }
\newcommand{\oud}{\<2_black> }
\newcommand{\out}{\<3_black>}
\newcommand{\outz}{\<K*3_black>}
\newcommand{\outu}{\<1K*3_black>}
\newcommand{\oudz}{\<K*2_black>}

\newcommand{\outd}{\<2K*3_black>}
\newcommand{\oudd}{\<2K*2_black>}
\newcommand{\oudX}{\<2x_black> }

\newcommand{\oub}{\<1_blue> }
\newcommand{\oudb}{\<2_blue> }
\newcommand{\outb}{\<3_blue>}
\newcommand{\outzb}{\<K*3_blue>}
\newcommand{\outub}{\<1K*3_blue>}
\newcommand{\oudzb}{\<K*2_blue>}

\newcommand{\outdb}{\<2K*3_blue>}
\newcommand{\ouddb}{\<2K*2_blue>}
\newcommand{\oudXb}{\<2x_blue> }

\newcommand{\oneb}{\blue{\mathbf{1}}}

\newcommand{\blue}[1]{\color{blue}#1\color{black}}

\newcommand{\lhs}{left hand side }
\newcommand{\rhs}{right hand side }
\newcommand{\hol}{Hölder }

\newcommand{\grad}{\triangledown}
\newcommand{\les}{\lesssim}
\renewcommand{\leq}{\leqslant}
\renewcommand{\geq}{\geqslant}
\renewcommand{\href}[1]{(\hyperref[#1]{\eqref*{#1}})}
\newcommand{\nsubset}{\not\subset}
\newcommand{\heat}{(\partial_t-\Delta)}

\newcommand{\R}{\mathbb{R}}
\newcommand{\E}{\mathbb{E}}

\newcommand{\eps}{\varepsilon}

\usepackage{color}
\numberwithin{equation}{section}
\newtheorem{thm}{Theorem}[section]
\newtheorem{cor}[thm]{Corollary}
\newtheorem{lem}[thm]{Lemma}

\theoremstyle{remark}

\newtheorem{rmq}[thm]{Remark}

\author{Augustin Moinat \and Hendrik Weber}

\address[Augustin Moinat]{University of Warwick, Coventry, United Kingdom}
\email{a.moinat@warwick.ac.uk}

\address[Hendrik Weber]{University of Bath, Bath, United Kingdom}
\email{h.weber@bath.ac.uk}

\thanks{HW is supported by the Royal Society through the University Research Fellowship UF140187.}
\thanks{The authors would like to thank the Isaac Newton Institute for Mathematical Sciences for its hospitality during the programme "Scaling limits, rough paths, quantum field theory" which was supported by EPSRC Grant Number: EP/R014604/1 }
\thanks{We would like to thank Ajay Chandra for useful feedback on a preliminary version of this article.}

\title{Space-time localisation for the dynamic $\Phi^4_3$ model }

\begin{document}

%\tableofcontents
%\newpage
\begin{abstract}
We prove an a priori bound for solutions of the dynamic $\Phi^4_3$ equation. 
This bound provides a control on solutions on a compact space-time set only in 
terms of the realisation of the noise on an enlargement of this set, and it does not
depend on any choice of space-time boundary conditions.

We treat the  large and small scale behaviour of solutions with completely different arguments. 
For small scales we use bounds akin to those presented in Hairer's theory of regularity 
structures. We stress immediately that our proof is fully self-contained, but we 
give a detailed explanation of how our arguments relate to Hairer's.
For large scales we use a PDE argument based on the maximum principle.
Both regimes are connected by a solution-dependent regularisation procedure.

The fact that our bounds do not depend on space-time boundary conditions makes them 
useful for the analysis of large scale properties of solutions. They can for example 
be used  in a compactness argument to construct solutions on the full space and their
invariant measures.

\bigskip



\end{abstract}
\maketitle

\section{Introduction}

The aim of this article is to derive  a priori bounds for the three dimensional stochastic quantisation equation, also known 
as the dynamic $\Phi^4_3$ model. This model is - at least formally - given by the non-linear stochastic partial differential 
equation (SPDE)
\begin{equation}\label{equation}
\heat u=-u^3+\xi,
\end{equation}
where $\xi$ is the space-time white noise over $\R \times \R^d$. In our main result, Theorem \ref{The Theorem}, we show a bound on the solution in the case $d=3$ on a compact space-time set that depends only on  a finite number of 
explicit polynomials in the Gaussian noise on a slightly larger space-time set. In particular, our bound does not depend on any space-time boundary conditions.

The main difficulty when working with \eqref{equation} is the roughness of the driving noise $\xi$ which in turn 
makes the solution irregular and the interpretation of non-linear terms non-trivial. It is now well-understood 
that solutions are distribution valued in spatial dimension $d \geq 2$  and the non-linearity has to be renormalised, 
which loosely speaking corresponds to replacing \eqref{equation} by 
\begin{equation}\label{renormalised}
\heat u=-u^3  +``\infty"\ u + \xi.
\end{equation}
The theory of singular SPDEs  of this type has been revolutionised in the recent years, starting with  Hairer's theory of regularity structures  \cite{hairer2014theory} 
and the theory of  paracontrolled distributions developed independently by Gubinelli, Imkeller and Perkowski \cite{Gubi}. 
Hairer's theory permits to develop a stable small scale theory, i.e. a local in time existence theory on compact spatial domains, 
 for a large class of SPDEs satisfying a scaling property called subcriticality. 
This notion corresponds exactly to super-renormalisability in quantum field theory and equation \eqref{equation} satisfies it for spatial dimension $d <4$.
In this approach, solutions are constructed in two steps. In the first step, a finite number of terms in a perturbative expansion of the solution based on the noise are 
constructed using probabilistic methods. In the second step, the actual solutions are sought in a space of distributions which are locally well approximated 
by these stochastic terms. The renormalisation procedure is treated in the probabilistic step making strong use of stochastic cancellations, while 
the second step is purely deterministic.
Hairer's work created a lot of activity. We mention in particular the works by Catellier and Chouk \cite{catellier2018}  and Kupiainen  \cite{Kupiainen2016} who produced similar 
short time existence and uniqueness results for \eqref{equation} using the method of paracontrolled distribution and  renormalisation group.

The theory of regularity structures is by now well-developed and permits to analyse a range of equations which are much more singular than the dynamic $\Phi^4_3$ model. It has 
been applied  to \eqref{equation} in "$4-\delta$" dimensions \cite[Section 2.8.2]{2017arXiv171110239B}, to the sine-Gordon model in the full subcritical regime \cite{2018arXiv180802594C}, 
and constructions of three dimensional gauge theories and the evolution of a random string on a manifold have been announced \cite{MartinFuture1},
\cite{MartinFuture2}.  
However, the arguments currently available are insufficient to go beyond a short time existence theory in any of these equations. 
For example, the construction of solutions to \eqref{equation} in \cite{hairer2014theory}  does not make use of the ``good'' sign of the nonlinear term $- u^3$ and would 
work equally if it were replaced with a $+u^3$. Solutions for this modified equation are expected to blow up in finite time.

 For \eqref{equation} in dimension $3$ the problem of passing from a local to a global solution theory has been largely overcome in a series of very recent works  starting with 
\cite{Mourrat3D} where \eqref{equation} was studied on the torus $\mathbb{T}^3$ and a priori estimates were obtained which ruled out the possibility of finite time blow-up.
In  \cite{GubinelliHofmanova1} a priori estimates for solutions on the full space $\R^3$ were shown;  see also \cite{albeverio2017invariant,GubinelliHofmanova2} for 
an analysis of the invariant measures based on similar ideas. All of these articles worked in the framework of  \emph{paracontrolled distribution} rather than regularity 
structures.

In this article we present a completely different technique to derive a priori estimates within the framework of regularity structures. We show a space-time version of the 
``coming down from infinity" property, i.e. we provide a bound on the solution on a compact space-time set that depends on the realisation of the noise on a slightly larger set, 
but does not depend on the behaviour of the solution elsewhere, making full use of the strong non-linear damping term $-u^3$. This local dependence makes this 
bound extremely useful when analysing the behaviour of solutions on large scales. 

A main interest of this approach is the technique itself. Its advantages are that we effectively separate the argument for small and large scales by dealing with a family of 
regularised equations for large scales and use (an appropriate restatement) of the theory of regularity structures
to analyse the small scales. This results in a relatively short argument compared to previous works and  
has the potential to work for a much larger class of equations.
We want to stress that our argument is fully self-contained and does not make use explicitly of any of the results in \cite{hairer2014theory}.
In fact, in both the statement of our main result and its proof we fully avoid the terminology of this theory, i.e. the  notions of model, modelled distribution, structure group etc.
but give a direct statement of all of the required bounds. This is possible, because the algebra involved in the small scale solution theory of \eqref{equation} 
is still not too complex and we hope that our direct approach  makes the presentation more clear. We do however include a separate section in which we translate our main estimates
into the regularity structure terminology.

% such as $\Phi^4_{4-\delta}$ (see \cite[Section 2.8.2]{2017arXiv171110239B})
%or potentially Gauge theories (see \cite{2018arXiv180104596S})

Finally, we would like to mention that in our companion paper  \cite{2018arXiv180810401M}  we have  implemented our approach in the case of one-dimensional reaction diffusion equations. 
Even in this much more regular case where no renormalisation enters, a priori bounds which do not depend on space-time boundary conditions seem to have been unknown.

This paper is structured as follows. Section~\ref{sec:Main} contains the elements needed to state our main result, Theorem~\ref{The Theorem}, starting with the definitions of the \hol right spaces in Section~\ref{mes_reg}.  Section~\ref{table tree} presents the setting in which we solve Equation~\eqref{equation} and our main result. The outline of the proof and the different lemmas required are presented in Section~\ref{outline}. We then explain the close connection between our setting and Hairer's theory of regularity structures in Section~\ref{connection}, where the full regularity structure for the $\Phi^4_3$ equation is presented. Section~\ref{The Proof} contains the proof of the Theorem \ref{The Theorem}, and Section~\ref{proofproof} contains the proof of the lemmas presented in Section~\ref{outline}.

\section{Setting and main result}\label{sec:Main}
\subsection{Measuring regularity}\label{mes_reg}

As usual when dealing with parabolic equations, regularity  will be measured with respect to the metric
\begin{equation}\label{parmetric}
d((t,x),(\bar{t}, \bar{x}))=\max\Big\{\sqrt{|t-\bar{t}|},|x-\bar{x}|\Big\},
\end{equation}
where $|\cdot|$ denotes the supremum norm on $\R^3$. 
We introduce the parabolic ball of center $z=(t,x)$ and radius $R$ in this metric $d$, looking only into the past: 
\begin{equation}\label{paraball}
B(z,R)=\{\bar{z}=(\bar{t},\bar{x})\in\R\times\R^3, d(z,\bar{z})<R,\bar{t}<t\}.
\end{equation}
We define the parabolic boundary of a set $D$ accordingly, as the set of points $z\in\bar{D}$ such that for any $r$, $B(z,r)\nsubset \bar{D}$. For $R>0$ we define $D_R\subset D$ as the set at distance $R$ from the parabolic boundary.
Set $P=(0,1)\times(-1,1)^3$ then we have
\begin{equation}\label{defpr}
P_R :=   (R^2,1) \times \{ x\colon |x|< (1-R) \}.
\end{equation} 
 Note that for $R'<R\leq 1$ we have for any domain $D$, $D_{R'}+B(0,R'-R)\subset D_R$. 
For $\alpha \in (0,1)$, we define the \hol semi-norm $[.]_\alpha$ 
%\comment{perhaps $z\prime$ instead of $y$?}
%
\begin{equation}\label{e:def-hol}
[u]_\alpha:=\sup_ {z\neq \bar{z}\in \R\times\R^3}\frac{|u(z)-u(\bar{z})|}{d(z,\bar{z})^\alpha}.
\end{equation}
For $\alpha \in (1,2)$, we define the \hol semi-norm $[.]_\alpha$  
\begin{equation}\label{e:def-hol2}
[u]_\alpha:=\sup_{\substack{z \neq \bar{z}  \in\R\times\R^3\\  z= (t,x); \; \bar{z} = (\bar{t}, \bar{z})}}\frac{|u(z)-u(\bar{z})-\grad u(z).(x-\bar{x})|}{d(z,\bar{z})^\alpha},
\end{equation}
where $\nabla$ refers to the spatial gradient.
 %\comment{not sure I like this convention. Why not call space time points $z=(t,x)$ and $z\prime = (t\prime, x\prime)$
%and then simply write $x-x\prime$?}
%\amcomment{At many points I use $x$ as a space-time base point}
%For a function $U$ of two variables, we have}
We will often deal with functions $U(z,\bar{z})$ of two variables 
generalising the increments of $u(z) - u(\bar{z})$ in \eqref{e:def-hol2} above. In this case we define for $\alpha \in (1,2)$

\begin{equation}\label{e:def-hol2var}
[U]_\alpha:=\sup_{\substack{z\in \R\times\R^3\\ z =(t,x)}}\inf_{\nu(z)\in \R^3}\sup_{\substack{\bar{z}\in \R\times\R^3\setminus\{z\} \\ \bar{z} = (\bar{t}, \bar{x}) }}\frac{|U(z,\bar{z})-\nu(z).(x-\bar{x})|}{d(z,\bar{z})^\alpha}.
\end{equation}
The infimum over functions $\nu$ is attained when  $\nu(z)$ is the spatial gradient in the second coordinate of $U$ at point $(z,z)$.
We often work with norms which only depend on the behaviour of functions / distributions on a fixed subset of time-space: if $B\subset\mathbb{R}\times\mathbb{R}^3$ is a bounded set, then we define the local $\alpha$-H\"older semi-norm  $[.]_{\alpha,B}$ as in \eqref{e:def-hol} with the supremum  restricted to $z,\bar{z}\in B$. The use of a third index $r$ as in $[.]_{\alpha,B,r}$ indicates that the supremum is restricted to $z$ and $\bar{z}$ at distance at most $r$.
Similarly, $\|.\|$ denotes the $L^\infty$ norm on the whole space $\R \times \R^3$ and $\|.\|_B$ the norm of the restriction of the function to $B$, and for a function of two variable, $\|.\|_{B,r}$ is the norm restricted to $z,\bar{z}\in B$ with $d(z,\bar{z})\leq r$.

 From now on  $x,y$ and $z$ will always denote a generic space-time variable, and we introduce the function $X$ which is the projection on space coordinates.

We work with a Besov-\hol type norm to measure negative regularity. These norms are usually defined by measuring the rate of blow-up when testing the distribution with a smooth approximation of Dirac. Different definitions usually boil down to different choices of approximations (e.g. convolution with rescaled smooth kernels in \cite[Definition 3.7]{hairer2014theory}, projection in Fourier space in Littlewood-Paley theory \cite{BCD} or wavelet basis
\cite{meyer1992wavelets}) and different choices usually require slightly different proofs of key properties such as multiplicative inequalities and Schauder
estimates. We need our testing operation to commute with the heat operator, which makes the convolution against rescaled kernels natural. 
Our definition is strongly inspired by the choice of smooth kernel satisfying the semi-group property with respect to the the scaling parameter
first introduced in \cite{OW}. This semi-group property allows to effectively connect regularisations at different scales and thus makes the proof of 
the Reconstruction Theorem \ref{Reconstruction} very convenient. However, an additional twist is required. 
For us it is important to be able to define local norms that only depend on properties of distributions on a compact set. 
This makes it most convenient 
to work with a compactly supported kernel in the definition of the norm.  But the kernel used in \cite{OW} does not have this property. Wavelet bases, on the other hand, 
permit a convenient transition from one scale to another and can consist of  compactly supported functions, but unfortunately the projection on these basis functions 
 do not commute with differential operators.  
The following simple 
construction yields a kernel which is compactly supported and enjoys a version of the semi-group property for dyadic scales which is enough 
to prove the reconstruction theorem.

We fix a non-negative smooth function $\Phi$ with support in $B(0,1)$, symmetric in space, 
%$\{(t,x)| (-t,x)\in B(0,1)\}$ 
with $\Phi(x)\in [0,1]$ for all $x\in \R\times\R^3$ and with integral $1$.
Setting $\Phi_T(t,x)=T^{-5}\Phi(\frac{t}{T^2},\frac{x}T)$,  
we now define $\Psi_{T,n}= \Phi_{T2^{-1}}\ast\Phi_{T2^{-2}}\ast...\ast\Phi_{T2^{-n}}$ and $\Psi_T=\lim_{n\rightarrow\infty}\Psi_{T,n}$ so that $\Psi_T=\Phi_\frac{T}2\ast\Psi_\frac{T}2$. The convergence can be checked easily. $\Psi_T$ and $\Psi_{T,n}$ are non-negative and smooth, symmetric in space and with support $B(0,1)$ and $B(0,1-2^n)$. We define the operator $(\cdot)_T$ by convolution with $\Psi_T$, and $(\cdot)_{T,n}$ by convolution with $\Psi_{T,n}$ for $n\geqslant 1$. $(\cdot)_{T,0}$ is the identity. Since $\Psi_{T,n+m}=\Psi_{T,n}\ast\Psi_{T2^{-n},m}$, we have
\begin{equation}\label{semigroup 2}
(\cdot)_{T,n+m}=((\cdot)_{T2^{-n},m})_{T,n}.
\end{equation}
Taking $m$ to infinity in this, or equivalently noticing that $\Psi_T=\Psi_{T,n}\ast\Psi_{T2^{-n}}$, we have the desired relation between dyadic scales
\begin{equation}\label{semigroup 1}
(\cdot)_T=((\cdot)_{T2^{-n}})_{T,n}.
\end{equation}
We then define the local $C^{\alpha}$ norm of a distribution $\theta$ for $\alpha<0$ as  
\begin{equation}\label{shauder ou}
[\theta]_{\alpha,C}=\sup_{T\leq 1}\|(\theta)_T\|_{C}T^{-\alpha}.
\end{equation}
It is proven in \cite[Theorem 2.34]{BCD} that for a similar quantity, in the case where $C$ is a torus of size one this corresponds to the classical Besov norm $\mathcal{B}^\alpha_{\infty,\infty}$ .
In our case, $[\theta]_{\alpha,C}$ depends on the distribution $\theta$ on $C+B(0,1)$ since $\Psi$ has support in $B(0,1)$

Furthermore, we mention the scaling estimates, for $n\in\mathbb{N}\cup\{\infty\}$ and $\alpha>-5$
\begin{equation}\label{moment of psi}
\int|\Psi_{T,n}(x-y)|d(x,y)^\alpha dy\leqslant T^\alpha,\quad \int|\grad\Psi_{T,n}(x-y)|d(x,y)^\alpha dy\les T^{\alpha-1}.
\end{equation}
Here and in the rest of the paper, "$\les$" denotes a bound that holds up to a multiplicative constant.
This immediately implies that for any $h\in C^\alpha$, $\alpha\in(0,2)$, and for any bounded set $C$,
we have  
\begin{equation}\label{mollified regularity 1}
\|h_T-h\|_{C_T}\leqslant T^\alpha\sup_{z\in C_T}[h]_{\alpha,B(z,T)}=T^\alpha[h]_{\alpha,C,2T}.
\end{equation}
Indeed, since $\Psi$ is symmetric in space we have $\int \Psi(y)X(y)dy=0$ and for all $x\in C$,
\begin{align*}
(h_T-h)(x)=&\int\Psi_T(x-y)(h(y)-h(x))dy\\
=&\int\Psi_T(x-y)(h(y)-h(x)-\ind{\alpha>1}\grad h(x).X(y-x))dy\\
\les& [h]_{\alpha,B(x,T)}\int\Psi_T(x-y)d(x,y)^\alpha dy.
\end{align*}
For products of function, we will sometimes be using the following notational convention:
\[
(fg)_T(x)-f(x)g_T(x)=((f-f(x))g)_T(x).
\]
The presence of the variable means that we evaluate the function there first, and the absence means that the convolution variable is used.

\subsection{Main result}\label{table tree}
We will work with a regularised version of \eqref{equation} throughout, i.e we assume that $u$ is a smooth function which on $P$ satisfies
\begin{equation}\label{phi43}
(\partial_t-\Delta)u=-u^3+\zeta+(3C_1-9C_2)u, 
\end{equation}
for real valued parameters $C_1, C_2$.  Thus, throughout the article, we never have to address the question how a given expression has to be interpreted to make sense. 
The main application we have in mind the case where  $\zeta = \xi_\delta$, i.e. a regularisation of the white noise at scale $\delta$ and where $C_1$ and $C_2$ 
are defined as the expectations of certain polynomials in $\xi_\delta$ which diverge like $\frac{1}{\delta}$ and $\log \delta^{-1}$ as the regularisation is removed. 
However, in our analysis these values only enter in the assumptions on the ``trees'' (see \eqref{def_tree2}, \eqref{def_oudd} and \eqref{def_outd}) and their
precise values do not appear. 
Despite dealing with smooth functions we stress that all of our estimates are stable in the limit $\delta \to 0$, where $\xi$ can only be measured as a distribution 
of regularity $- \frac52-$ and $u$ as a distribution of regularity $-\frac12-$. We will freely use the convention to speak of ``distributions'' when we refer to smooth functions
that can only be measured in a distributional norm in this limit.

%
%
%In order to prove a theorem that holds for solutions to the equation , we actually consider a smooth version of it, and prove bounds that hold independently of the space-time regularisation, thus yielding bounds that also hold in the limit where that small scale cut-off is removed. We denote by $\zeta$ a regularisation of the white noise: $\zeta=\xi_\delta$ for some $\delta>0$. The renormalised equation \eqref{renormalised} is then restated as
%\begin{equation}\label{phi43}
%(\partial_t-\Delta)u=-u^3+\zeta+Cu, 
%\end{equation}
%where $C$ is a constant that diverges to infinity when $\delta$ goes to zero. We assume that this holds point-wise in the box $P$. The regularisation parameter $\delta$ will not be made explicit in the future, but all bounds we prove hold uniformly in $\delta$.
%%

We first introduce several polynomials in $\zeta$ that are used in the local description of the solution to \eqref{phi43}. 
These are (essentially) the same objects which appear in Hairer's small scale solution theory for \eqref{phi43} and we use his convention 
to denote these objects by trees.

%
% functions and distributions used for a local description of the solution to \eqref{phi43}, represented by trees.  A dot denotes one occurrence of the regularised noise $\zeta$, and a vertical line denotes integration against the heat kernel. Joining two trees at the root denotes a renormalised multiplication. Subcriticality tells us that only a finite number of such trees should be enough to describe our solution to a high enough regularity.

%
%These are the same (up to a tiny technicality) trees encountered in Hairer's theory of regularity structure. 
%Bounds on these objects have been derived in \cite[Section 10]{hairer2014theory} (and for closely 
%related objects also in \cite{catellier2018,Pedestrians}). Here we treat the bounds on these objects as input to our theory.
%%

We start by fixing an $\eps>0$ which will always be assumed to be ``sufficiently small".
The first tree is $\ou$ is assumed to  satisfy the point-wise identity on  $P$
\begin{equation}\label{def_tree1}
\heat\ou=\zeta,
\end{equation}
and we assume a control in the   $C^{-\frac12-\epsilon}$ norm. 
%
% Subtracting this linear part however yields an object expected to be of higher regularity $v:=u-\ou\in C^{\frac12-3\epsilon}$ uniformly in the regularisation as will be shown below.
%
%We decompose the constant $C$ into $C:=3C_1-9C_2$. 
The constant $C_1$ appears in the following definitions:
\begin{equation}\label{def_tree2}
\oud:=\ou^2-C_1,\quad\out:=\ou^3-3C_1\ou,
\end{equation}
and these distributions will be measured in the norms of $C^{-1-2\eps}$ and $C^{-\frac32-3\eps}$, respectively.
We also introduce symbols of higher order: We assume that $\oudz$ and $\outz$ satisfy the point-wise identity on $P$
\begin{equation}\label{def_tree3}
\heat\oudz=\oud,\quad\heat\outz=\out.
\end{equation}
As expected with the heat operator, we assume that the regularity  is increased by $2$ i.e. $\oudz\in C^{1-2\epsilon}$ and $\outz\in C^{\frac12-3\epsilon}$.
Finally we introduce the trees $\oudX$ denoting the product of $\oud$ with $X$, $\oudd,\outu$ and $\outd$ and for these we will need  assume bounds on the quantitites 
\begin{align}
[\oudX]_{-2\epsilon}&=\sup_{x\in P}\sup_{T<1}T^{2\epsilon}\Big|\int X(y-x)\oud(y)\Psi_T(y-x)dy\Big|,\label{def_oudX}\\
[\oudd]_{-4\epsilon}&=\sup_{x\in P}\sup_{T<1}T^{4\epsilon}\Big|\int(\oudz(y)\oud(y)-C_2-\oudz(x)\oud(y))\Psi_T(y-x)dy\Big|,\label{def_oudd}\\
[\outu]_{-4\epsilon}&=\sup_{x\in P}\sup_{T<1}T^{4\epsilon}\Big|\int(\outz(y)\ou(y)-\outz(x)\ou(y))\Psi_T(y-x)dy\Big|,\label{def_outu}\\
[\outd]_{-\frac12-5\epsilon}&=\sup_{x\in P}\sup_{T<1}T^{\frac12+5\epsilon}\Big|\int(\outz(y)\oud(y)-3C_2\ou(y)-\outz(x)\oud(y))\Psi_T(y-x)dy\Big|\label{def_outd}.
\end{align}

We will work with the function  $v:=u-\ou$ which satisfies
\begin{align}\label{eq v}
\heat v=&-u^3+(3C_1-9C_2)u=-(v+\ou)^3+(3C_1-9C_2)(v+\ou)\nonumber\\
=&-v^3-3v^2\ou-3v(\ou^2-C_1)-(\ou^3-3C_1\ou)-9C_2(v+\ou)\nonumber\\
=&-v^3-3v^2\ou-3v\oud-\out-9C_2(v+\ou).
\end{align}
The fact that the constant $C_1$ disappears in this expansion was already noted in \cite{DPD} and that was enough to define solutions in dimension $2$, where the constant $C_2$ is unnecessary.

%
%In terms of regularity the power counting can be summarized as such: the white noise is of regularity $--\frac{d+2}2-\epsilon=\frac52-\epsilon$; multiplication adds regularities and integration against the heat kernel adds $2$ regularity. 
%This corresponds to the notion of homogeneity in Hairer's theory of regularity structures, and the decomposition in \eqref{def_oudX}-\eqref{def_outd} corresponds to positive renormalisation encoded in Hairer's work in the operator  $\Pi_x$.

\begin{thm}\label{The Theorem}
If $v$ solves \eqref{eq v} pointwise on $P$, then we have:  
\begin{equation}\label{equation of the theorem}
\|v\|_{P_R}\leqslant C\max\Big{\{}\frac1{R},[\tau]_{|\tau|}^\frac1{n_\tau(\frac12-\epsilon)},{\tau\in L}\Big{\}},
\end{equation}
where $L=\{ \ou ,\oud ,\oudX, \oudz, \out, \outz, \oudd , \outu, \outd  \}$ and $|\tau|$ is the regularity in which we measure 
the tree $\tau$ in the way explained above.
\end{thm}
%The proof of this theorem relies on the following lemmas.

%
\begin{rmq}

As stated above, the bounds we assume on the ``trees'' are (almost -- see the following Remarks~\ref{The Remark} and ~\ref{The Remark 2})
identical to those appearing as  input into the analytic part of \cite{hairer2014theory}. The particular form of 
the $x$-dependent ``counterterms'' $-X(x)\oud(y)$ in \eqref{def_oudX}, $ \oudz(x)\oud(y)$ in \eqref{def_oudd}
and $-\outz(x)\oud(y)$ in \eqref{def_outu} and $\outz(x)\oud(y)$ in \eqref{def_outd} corresponds 
exactly to the ``positive renormalisation'' or re-centering procedure of the trees performed there.
%The constants $C_1$ and $C_2$ appearing in the definition correspond to the ``negative renormalisation''
%that removes divergencies. 
See Section~\ref{sec:RG} for a more detailed discussion of positive %and negative 
renormalisation in the theory of regularity structures.

In the case where $\zeta = \xi_\delta$ is a regularised white noise and where $C_1 = \E \ou (y)^2$ and $C_2 = \E \oudz(y)\oud(y)$,
e.g. for $y = (1,0)$, uniform-in-$\delta$-bounds on the various norms were obtained in \cite[Section 10]{hairer2014theory}.
We stress that in this low-regularity situation the convergence of these terms as $\delta \to 0$ is highly non-obvious, even after 
renormalisation. The calculations use probabilistic tools and strongly rely on stochastic cancellations. 

%
%Our argument takes all the definitions of the trees as input for some given value of $C_1$ and $C_2$. 
%When we go back to the original equation with white noise, 
%the convergence of the renormalised products in the given space is far from obvious due to their low regularity and the argument in \cite{catellier2018,Pedestrians} crucially exploits stochastic cancellations. 

The estimates in \cite[Section 10]{hairer2014theory} actually yield bounds on the moments of all of  these terms, so that our main result
\eqref{equation of the theorem} implies bounds  moments of the solution: for each  $\tau$ is $[\tau]_{\gamma_\tau}$ is a random
variable in the (inhomogeneous) Wiener chaos of order $n_\tau$ over the Gaussian noise, so that one gets 
\[
\E[\exp (\lambda[\tau]_{\gamma_\tau}^{\frac{2}{n_\tau}})]<\infty,
\]
for some $\lambda>0$. Hence for $\bar{\lambda} =\frac{\lambda}{C^{1-2\epsilon}} $ we get
\begin{equation}\label{integrabilitybound}
\E[\exp (\bar{\lambda} \|v\|_{P_R}^{1-2\epsilon} )]<\infty.
\end{equation}

\end{rmq}

\begin{rmq}
One of the main motivations to consider  \eqref{equation} is to use the Markovian dynamics described by it to study its invariant measure,
the Euclidean $\Phi^4_3$ quantum field theory. In order to link this Euclidean (imaginary time) field theory to a real time field 
theory, this measure should satisfy certain properties, the Osterwalder-Schrader axioms \cite[section 6.1]{glimm2012quantum}.
Our bounds \eqref{integrabilitybound} immediately transfers to this invariant measure. Unfortunately, these stretched exponential 
moments just fall short of  the exponential bounds required for the Analyticity Axiom.

\end{rmq}

%\begin{rmq}\label{rem2}
%A first result that follows directly from Theorem \ref{The Theorem} is the stochastic integrability of the solution in the white noise limit. When constructing the trees one can obtain bounds on the $p$-th moments of these, as is done in \cite[Section 10]{hairer2014theory} using a Wiener chaos decomposition. One can then track the value of the constant to obtain exponential integrability is related to $n_\tau$:
%\[
%\E[\exp (\lambda[\tau]_{\gamma_\tau}^{\frac{2}{n_\tau}})]<\infty,
%\]
%for some $\lambda>0$. This is done in \cite[Section 3]{Pedestrians}. Hence
%\[
%\E[\exp (\lambda\|v\|_{P_R})^{1-2\epsilon}]<\infty.
%\]
%Unfortunately this falls just short of exponential bounds required for the Analyticity in Osterwalder-Schrader axioms \cite[section 6.1]{glimm2012quantum}.
%\end{rmq}

\begin{rmq}\label{The Remark}

Hairer's convention in the definition of the symbols $\ou$ in \eqref{def_tree1}, $\oudz$ and $\outz$ in \eqref{def_tree3} differs slightly from ours: Instead of assuming 
that these objects satisfy a partial differential equation as we do, he defines them using an integral condition, e.g.
\[
\ou(x) = \int_{\R \times \R^3} K(x,y) \zeta(y) dy,
\]
for a singular integral kernel $K$. This kernel $K$ is essentially the Gaussian heat kernel, but it is post-processed to make it compactly supported and to 
integrate to $0$ against  polynomials up to a certain degree. After this post-processing $K$ is not associated to a differential operator any more and in this 
definition $\ou$ and the other stochastic terms are not characterised by a (simple) PDE.  This is in line with the general philosophy pursued in \cite{hairer2014theory}
 to view \eqref{equation} as an integral equation using the mild formulation rather than a differential equation.
 
\end{rmq}

\begin{rmq}\label{The Remark 2}
Continuing the discussion of the symbols $\ou$, $\oudz$ and  $\outz$  we point out that in \eqref{def_tree1} and  \eqref{def_tree3} we do not impose boundary conditions, 
but only that a certain PDE holds point-wise. There is thus some choice in how these objects are defined and our main result, the estimate \eqref{equation of the theorem}, holds uniformly 
over all of these choices. This is also the reason why the symbols $\oudz$ and $\outz$ appear in the list $L$. For many choices of boundary conditions 
Schauder theory would imply $[\oudz]_{1-2\epsilon} \lesssim [\oud]_{-1-2\epsilon}$ and $[\outz]_{\frac12-3\epsilon} \lesssim [\out]_{-\frac32-3\epsilon}$ so that these symbols could be removed from $L$. 

A natural choice to would be to impose Dirichlet boundary conditions on the parabolic boundary of $P$ in \eqref{def_tree1} and \eqref{def_tree3} and in this case such a Schauder estimate holds indeed. 
Moreover, with this choice one would have the nice property that all of the objects on the right hand side only depend on the realisation of $\zeta$ on $P$, which would be in line with a ``space-time Markov
property.'' This nice choice has the slight disadvantage that (in the case where $\zeta = \xi_\delta$ is a regularised white noise) the negative renormalisation would have to be 
modified reflecting the boundary conditions which would lead to $x$-dependent $C_1$ and $C_2$ in \eqref{def_tree2}, \eqref{def_oudd}, \eqref{def_outd}, and then an extra term would have to be added 
in \eqref{phi43} in order to make the renormalisation of the original equation $x$-independent. Such a construction could certainly be implemented, but we refrain from doing so here (see \cite{gerencser2017singular},
\cite{TW} for discussion of similar boundary issues).

\end{rmq}

%
%
%The definition given of the integrated stochastic objects $\ou, \outz$ and $\oudz$ was not precise in terms of boundary conditions, although we mentioned estimates on the norms, namely that $\ou\in C^{-\frac12-\epsilon},\oudz\in C^{1-2\epsilon}$ and $\outz\in C^{\frac12-3\epsilon}$. If one chooses for example Dirichlet boundary conditions on the set $P$, we have from a Schauder estimate $[\oudz]_{1-2\epsilon}\les[\oud]_{-1-2\epsilon}$ and $[\outz]_{1-2\epsilon}\les[\out]_{-1-2\epsilon}$
%(see e.g. \cite[Lemma 5.1]{2018arXiv180810401M}). However in that case the renormalisation procedure to define the products requires the constant $C_1$ and $C_2$ to be replaced by space-time dependent functions, vanishing near the boundary. Allowing for this would have introduced some technicalities which would hinder the clarity of the proof. A cheap way of keeping $C_1$ and $C_2$ as constants is to assume Dirichlet boundary conditions on a bigger set, for example $P+B(0,1)$. In that case, one looses the good Markov property that follows from the strong localisation of the theorem and gets a bound on $\|v\|_{P_R}$ depending on the noise on the set $P+B(0,1)$.
%%\end{rmq}
\begin{rmq}
The spatial dimension $d=3$ only enters our analysis through the regularity assumptions on the ``trees''. The various $\gamma_{\tau}$ are all derived from 
the parabolic regularity of the white noise in $1+3$ time-space dimensions , which is $- \frac{5}{2}-$. The actual PDE arguments we present do not rely on a specific 
choice of $d$.

%
%The dimension only enters in this result through the regularity assumption on the white noise. In particular in the proof only the values of some constants would be changed in higher dimension, if we assumed that we had a driving noise of the same regularity which allowed to prove convergence of the desired products.
\end{rmq}

\subsection{Outline of proof}\label{outline}

One of the key ideas behind the theory of regularity structures is the following scaling argument:
\begin{equation}\label{SHE scaling}
\hat{u}(t,x) = \lambda^{\frac{d}{2}}u(\lambda^2 t , \lambda x)
\end{equation} 
is the scaling under which the stochastic heat equation $\heat u=\xi$ is invariant in law. For the $\Phi^4$ equation the non-linearity $-u^3$ scales like $\lambda^{4-d}\hat{u}^3$. In dimension less than $4$, this term formally vanishes 
on small scales, i.e. when $\lambda$ goes to zero. This property is called subcriticality in Hairer's theory and corresponds to super-renormalisability in quantum field theory. 
This observation suggests that in order to control the behaviour of $u$ on ``small scales'' one should use the heat operator and treat the non-linearity as a perturbation. This
 is precisely how a small-scale local solution theory is built in \cite{hairer2014theory}. The sign of the non-linearity $-u^3$ is not used in this argument. The argument for large scales on the other side 
 clearly has to rely on the ``good term'' $-u^3$ and should not use the smoothing of the heat operator too much.

We have already seen that as a perturbation of the linear equation, $v=u-\ou$ satisfies
\begin{equation}\label{eq v bis}
\heat v=-v^3-3v^2\ou-3v\oud-\out-9C_2(v+\ou).
\end{equation} 
To control large scales, we apply the regularising operator $( \cdot)_T$ for some $T$ to be chosen below and we get the equation
\begin{equation}\label{phi43_T}
(\partial_t-\Delta)v_T=-(v_T)^3-3(v^2\ou)_T-3(v\oud)_T-9C_2(v_T+\ou_T)-(\out)_T+((v_T)^3-(v^3)_T).
\end{equation}
This equation is not closed in terms of $v_T$ and we will require control on the commutator $(v_T)^3-(v^3)_T$ and on the products $(v^2\ou)_T$ and $(v\oud)_T$. These are bounded in the small-scale theory. For large scale bounds, we use the following lemma
\begin{lem}\label{main thm poly 1}
Let $u$ be a continuous function defined on $[0,1]\times [-1,1]^3$, for which the following holds point-wise in $(0,1]\times (-1,1)^3$:
\begin{equation}\label{eq:max rd 1}
(\partial_t-\Delta)u=-u^3+g(u,z),
\end{equation}
where $g$ is a bounded function.
We have the following point-wise bound on $u$, for all $(t,x)\in(0,1]\times (-1,1)^3$:
\begin{equation}\label{eq:max theorem 1}
| u(x,t) | \leqslant C\max\Big\{\frac1{\min\{\sqrt{t},(1-x_i),(1+x_i),i=1,2,3\}},\|g\|^\frac13\Big\},
\end{equation}
for some independent constant $C$.
\end{lem}
This lemma is a simplified version of \cite[Theorem 4.4]{2018arXiv180810401M} and the proof (in Section \ref{proof max}) is based on the maximum principle. It is the only part of the argument which makes use of the fact that $u$ is a scalar field and not vector valued. The rest of 
the proof would go through in the vector-valued case and we expect that it is possible to find a vector-valued replacement for Lemma 2.6 as well.

In order to close the estimate obtained from Lemma \ref{main thm poly 1}, we require a bound that allows to 
control high order regularity of $v$ in terms of the $L^\infty$ norm. The classical method would consist of using Schauder estimate of the form \cite[Theorem 8.9.2]{krylov1996lectures}
\[
[u]_{\delta+2,D_R}R^{2} \lesssim [\heat u+u]_{\delta,D} 
\]
for solutions of inhomogeneous heat equation.  Then if the \rhs depends on a lower order norms of $u$ it can be absorbed into the \lhs.
%
%This is close in spirit to the "fixed point problem" of local solution theory. 
We perform such an argument in the case where usual \hol norms are replaced by the norms of "modelled distributions" 
(which depend on the underlying noise $\zeta$).

First,  power counting suggests that $v+ \outz$ has better regularity than $v$ (namely $1-2\epsilon$)  and that this would be enough 
to define $v^2 \ou=v(v+\outz)\ou-\outz\ou$ (assuming that we can construct $\outz \ou$), but not enough to define $v \oud$. 
The next idea to get even better description of solution by explicit stochastic terms is to freeze coefficients at base point, and to look at local expansions that depend on that base point. The expansion of $v$ in around base point $x$ goes as follows:
\begin{equation}\label{eq:local-description-v}
v(y)=v(x)-(\outz(y)-\outz(x))-3v(x)(\oudz(y)-\oudz(x)).
\end{equation}
We introduce the following function of two variables based on this local description:
\begin{equation}\label{defU}
U(x,y)=v(y)-v(x)+\outz(y)-\outz(x)+3v(x)(\oudz(y)-\oudz(x)).
\end{equation}
The regularity of $U$, as defined in \eqref{e:def-hol2var}, is expected be higher than $1$. This better description is indeed enough to define $v \oud$. The core observation is the following abstract 
reconstruction theorem, which is a variant of  \cite[Theorem 3.10]{hairer2014theory} and \cite[Proposition 1]{otto2018parabolic}. 
%Since we are only dealing with functions as we have regularised the equation, we do not concern ourselves with any existence statement.
\begin{thm}[Reconstruction]\label{Reconstruction}
Let $\gamma>0$ and $A$ be a finite subset of $(-\infty,\gamma]$. Let $T\in (0,1)$ and $x\in \R\times \R^3$. 
For a function $F \colon B(x,T)^2 \to \R$ assume that for all $\beta\in A$ there exist constants $C_\beta>0$ and $\gamma_\beta\geqslant\gamma$ such that for all $t\in(0,T)$, for all $x_1,x_2\in B(x,T-t)$
\begin{equation}\label{eq:reconstruction_hypothesis}
\Big| \int\Psi_t(x_2-y)(F(x_1,y)-F(x_2,y))dy \Big|\leq \sum_{\beta\in A}C_\beta d(x_1,x_2)^{\gamma_\beta-\beta}t^\beta.
\end{equation}
Then $f:y\mapsto F(y,y)$ satisfies
\begin{equation}\label{eq:reconstruction_conclusion}
\Big| \int\Psi_T(x-y)(F(x,y)-f(y))dy \Big|\les\sum_{\beta\in A}C_\beta T^{\gamma_\beta},
\end{equation}
where ``$\les$'' represents a bound up to a multiplicative constant depending only on $\gamma$ and $A$.
\end{thm}
As a consequence of this theorem, we get the following bounds on the products.
\begin{lem}\label{lem_v2ou}
The following bound on $v^2\ou$ holds:
\begin{align}\label{eq:lem_v2ou}
|(v^2\ou)_T(x)|\les& T^{\frac12-3\epsilon}\|v\|_{B(x,T)}[v+\outz]_{1-2\epsilon,B(x,T) } [\ou]_{-\frac12-\epsilon}\nonumber\\
&+T^{\frac12-7\epsilon}[v]_{\frac12-3\epsilon,B(x,T)}([\outz]_{\frac12-3\epsilon}[\ou]_{-\frac12-\epsilon}+[\outu]_{-4\epsilon})\\
&+\|v\|_{B(x,T)}^2[\ou]_{-\frac12-\epsilon}T^{-\frac12-\epsilon}+ \|v\|_{B(x,T)}[\outu]_{-4\epsilon}T^{-4\epsilon}\nonumber.
\end{align}
\end{lem}
\begin{lem}\label{lem_voud}
The following bound on $v\oud$ holds:
\begin{align}\label{eq:lem_voud}
|((v-v&(x))\oud)_T(x)+3C_2(v_T+\ou_T)(x)|\les \nonumber\\
T^{\frac12-7\epsilon}\Big(&[v]_{\frac12-3\epsilon,B(x,T)}[\oudd]_{-4\epsilon}
+[U]_{\frac32-5\epsilon,B(x,T)} [\oud]_{-1-2\epsilon}+[\nu]_{\frac12-5\epsilon,B(x,T)}[\oudX]_{-2\epsilon}\Big)\nonumber\\
&+[\outd]_{-\frac12-5\epsilon}T^{-\frac12-5\epsilon}+T^{-4\epsilon}\|v\|_{B(x,T)}[\oudd]_{-4\epsilon}+\|\nu\|_{B(x,T)}[\oudX]_{-2\epsilon}T^{-2\epsilon},
\end{align}
where $U$ is as introduced in \eqref{defU} and $\nu$ is optimal in the definition \ref{e:def-hol2var}.
\end{lem}
The proof of Theorem~\ref{Reconstruction} and the two Lemmas \ref{lem_v2ou} and \ref{lem_voud} can be found  in Section \ref{proof reconstruction}.
To bound the quantities appearing in the \rhs of these lemmas, we introduce our version of the Schauder estimate.
\begin{lem}\label{lemschauder}
Let $1<\kappa<2$ and $A\subset(-\infty,\kappa]$ finite. Let
$U$ be a bounded function of two variables defined on a domain $D\times D$ such that $U(x,x)=0$ for all $x$. Let $d_0>0$ and assume that for any $0<d\leq d_0$ and $L\leq \frac{d}4$ there exists a constant $M^{(1)}_{D_d,L}$ such that for all base points $x\in D_d$ and length scales $T\leqslant L$, it holds that
\begin{equation}\label{lemschauder1}
T^2\|\heat U_T(x,\cdot)\|_{B(x,L)}\leqslant M^{(1)}_{D_d,L}\sum_{\beta\in A}T^\beta L^{\kappa-\beta}.
\end{equation}
Assume furthermore, that for $L_1,L_2\leq \frac{d}4$ there exists a constant $M^{(2)}_{D_d,L_1,L_2}$ such that
for any $x\in D_d$, for any $y\in B(x,L_1)$, for any $z\in B(y,L_2)$ the following "three-point continuity" holds:
\begin{equation}\label{lemschauder2}
|U(x,z)-U(x,y)-U(y,z)|\leqslant M^{(2)}_{D_d,L_1,L_2}\sum_{\beta\in A}d(y,x)^\beta d(z, y)^{\kappa-\beta}.
\end{equation}
Additionally define
\[
M^{(1)}:=\sup_{d\leq d_0}d^\kappa M^{(1)}_{D_d,\frac{d}2},\quad\text{and}\quad M^{(2)}:=\sup_{d\leq d_0}d^\kappa M^{(2)}_{D_d,\frac{d}2,\frac{d}4}.
\]
Then
\begin{equation}\label{lemschauderC}
\sup_{d\leq d_0}d^\kappa [U]_{\kappa,D_d}\les M^{(1)}+M^{(2)} +\sup_{d\leq d_0}\|U\|_{D_d,d}.
\end{equation}
Here and in the proof,"$\les$" denotes a bound that holds up to a multiplicative constant that only depends on $\kappa$ and $A$.
\end{lem}

\begin{cor}\label{corschauder}
Fix $0<d\leq d_0$ such that $D_d\neq\emptyset$. Assume that $D_d$ satisfies a spatial interior cone condition with parameters $r_d>0$ and $\lambda\in(0,1)$, i.e. for all $r\in[0,r_d]$, for all $x\in D_d$, for any vector $\nu\in \R^3$, there exists $y\in D_d$ such that $d(x,y)=r$
\[
|\nu.X(y-x)|\geq \lambda |\nu|d(x,y).
\]
Then for a (near) optimal function $\nu$ in \eqref{lemschauderC}, for all $r\in[0,r_d]$,
\begin{equation}\label{corschauder1}
\lambda\|\nu\|_{D_d}\leq [U]_{\kappa,D_d}r^{\kappa-1}+\|U\|_{D_d,r}r^{-1}.
\end{equation}
If \eqref{lemschauder2} holds for all $x,y,z\in D_d$ we have for $r\leq r_d$, 
\begin{equation}\label{corschauder2}
[\nu]_{\kappa-1,D_d}\lesssim [U]_{\kappa,D_d}+M^{(2)}_{D_d,\frac{d}4,\frac{d}4}+r^{-\kappa}\|U\|_{D_d,r}.
\end{equation}
Here and in the proof,"$\les$" denotes a bound that holds up to a multiplicative constant that only depends on $\kappa$.
\end{cor}
Note that if $D=P_R$ for $R<\frac12$, the interior cone condition holds for $D_d$ with $r_d= \frac12-d$ and $\lambda=\frac{\sqrt{2}}2$, uniformly in $R$. The proof of the lemma and its corollary are in Section \ref{proof schauder}

The Schauder estimate and the reconstruction lemmas can be combined with the large scale bound into a self-consistent bound that can be iterated leading to our main result, Theorem~\ref{The Theorem}. 
This argument can be found in Section~\ref{The Proof}.

\section{Translation to the language of regularity structures}\label{connection}
\label{sec:RG}

Although our argument is not formulated using the terminology of the theory of regularity structures, the analysis of the  small scale behaviour,  Theorem~\ref{Reconstruction},
Lemmas~\ref{lem_v2ou} and \ref{lem_voud} as well as the Schauder estimated~\ref{lemschauder}, builds on the same key ideas as this theory.
We now provide a translation of how the lemmas that appear in our paper can be stated in terms of the central objects introduced in the theory of regularity structures such as the models, modelled distributions, and the abstract integration operator.

We begin by recalling the setup in \cite{hairer2014theory}:
In \cite[Definition 2.1]{hairer2014theory} a regularity structure is defined as a triple $(A,T,G)$, consisting of an index set $A \subset \R$, a graded vector space $T = \bigoplus_{\alpha \in A} T_\alpha$ and a group
$G$ of linear  transformations acting on $T$ with some additional properties. 
In this framework, the local description of the solution $u$ is encoded by replacing the scalar valued function / distribution $u$ by modelled distribution, which is a function $\mathcal{U}:\R\times\R^3\to T$ for a certain purpose-built regularity structure. To build this structure one first introduces some  symbols, namely
\[
\{ \oneb, \oub, \outzb, \oudzb \} \cup \{ \blue{X_i} : i=1,2,3 \} .
\]
At this level these blue symbols are completely abstract objects, but of course they ultimately represent the functions / distributions appearing in the local description.
To each of these symbols $\tau$ one associates a \emph{homogeneity} $|\tau| \in \R$, namely 
\[
| \oneb | = 0, \quad | \blue{X_i}| = 1, \quad  | \oub | = -\frac{1}{2}-\epsilon, \quad | \outzb| = \frac12 - 3\epsilon, \quad |\oudzb| = 1-2\epsilon.
\]
The space   $T$ is then defined as the finite dimensional space 
\[
T = \bigoplus_{\tau \in \{ \oneb, \blue{X_i}, \oub, \outzb, \oudzb \} }  \R \tau,
\]
and $A$ is defined to be the set of homogeneities of these symbols. It turns out that the modelled distribution $\mathcal{U}$ takes the form 
\begin{equation} \label{def_phi_regstr}
\mathcal{U}(x) = \oub + v(x) \oneb - \outzb  - 3v(x) \oudzb- \nu(x) . \blue{X},
\end{equation}
for some functions  $v$ and $\nu$ (which of course  coincide with our functions $v$ and $\nu$). 
For our analysis we choose to work with a local description for $v = u - \ou$,
 which in the notation of regularity structures would take the form
\begin{equation} \label{def_V_regstr}
\mathcal{V}(x) =  v(x) \oneb - \outzb  - 3v(x) \oudzb- \nu(x) . \blue{X},
\end{equation}
i.e. the only difference with respect to \eqref{def_phi_regstr} is that the term $\oub$ is removed. 
Equation \eqref{def_V_regstr} should be viewed as an abstract counterpart of our equation~\eqref{eq:local-description-v}.
For us it is more convenient to work with $v$ rather than $u$ to get good bounds on the error term 
$(v_T)^3- (v^3)_T$. We argue below, that the regularity assumption we impose on $\mathcal{V}$ is equivalent 
to the condition imposed on $\mathcal{U}$ in \cite{hairer2014theory}.

Just like our main  result, Theorem~\ref{The Theorem}, the solution theory using regularity structures requires a perturbative expansion as an input. There this
expansion is encoded in the notion of a model \cite[Definition 2.17]{hairer2014theory}. To each of the symbols, one associates a function / distribution $\mathbf{\Pi}\tau$ corresponding exactly to our definitions \eqref{def_tree1},
 and \eqref{def_tree3}, i.e.
 \begin{align*}
 \mathbf{\Pi} \oneb (y) &= 1, &
\mathbf{\Pi} \blue{X_i}(y) &= y_i, &    
 \mathbf{\Pi} \oub (y) &= \ou(y), \\
 \mathbf{\Pi} \outzb (y) &= \outz(y), &
 \mathbf{\Pi} \oudzb (y) &= \oudz(y). &
 \end{align*}
A key idea of the theory is  to not work with these  distributions directly, but with \emph{centered} or \emph{positively renormalised} objects, $\Pi_x \tau$
indexed by a base-point in $x \in \R \times \R^3$. The right notion of regularity for the  \emph{modelled distributions} $\mathcal{U}$ and $\mathcal{V}$ is then defined in term of this recentering procedure.  
%Later in this section we show how the group $G$ encodes this recentering.  \comment{reformulate}

For the symbols we have introduced so far, the centering is relatively simple  and amounts to subtracting  the value at the base point for the symbols of strictly positive homogeneity:
\begin{align*}
\Pi_x\oneb(y)&=1, & \Pi_x \blue{X_i}(y) &= y_i-x_i , &
\Pi_x \oub (y) &= \ou(y) ,&\\
\Pi_x \oudzb(y)&=\oudz(y)-\oudz(x), &
\Pi_x\outzb(y)&=\outz(y)-\outz(x).
\end{align*}
The reason why one works with these centered objects is that one has good control over their behaviour as the argument approaches the base point $x$. This is encoded 
in the formula \cite[Equation (2.15)]{hairer2014theory},
\begin{equation}\label{model-continuity}
\langle \Pi_x \tau, \varphi^\lambda_x \rangle \lesssim \lambda^{|\tau|},
\end{equation}
where $\varphi$ is a smooth test-function rescaled to scale $\lambda$ and centered at the base-point $x$. This corresponds exactly to our 
regularity assumption on the H\"older norms of the objects, see Section~\ref{table tree} (our scale is called $T$ rather than $\lambda$ and the test-function is called $\Psi$ rather than $\varphi$).
%We stress right away that the centering of the additional symbols appearing in the definition of products below does not reduce to subtracting the value at the base point and that for these symbols condition 
%\eqref{model-continuity} becomes strictly stronger than a simple H\"older continuity statement.
%
% 
%
%

In order to connect the centering procedure to the functions  $\mathcal{U}$ and $\mathcal{V}$  and to formulate  the right continuity condition,  it is useful to introduce the structure group $G$.
In the current context this group is simply the five-dimensional group of  all linear transformations $F$ on $T$ of the form
\begin{align} 
\notag
F \oneb &= \oneb, &
F \blue{X_i} &= \blue{X_i} + a_i \oneb  \quad a_i \in \R, &
F \oub & = \oub, \\
\label{def_group}
F \outzb &= \outzb + b \oneb \quad  b \in \R,  &
 F \oudzb &= \oudzb + c \oneb  \quad c \in \R,
\end{align}
but this group will be enlarged as more symbols are introduced below.
For each $x \in \R \times \R^3 $ we define $F_x \in G$ by 
\begin{align*}
F_x\oneb&= \oneb , & F_x \blue{X_i}(y) &= \blue{X_i} - x_i \oneb  , &
F_x \oub  &= \oub  ,&\\
F_x \oudzb&=\oudzb -\oudz(x) \oneb, &
F_x\outzb&=\outzb -\outz(x) \oneb,
\end{align*}
so that one  gets 
\[
\Pi_x \tau = \mathbf{\Pi} F_x \tau. 
\]
Now, for $x,y \in \R \times \R^3$ we set 
\[
\Gamma_{xy} = F_y^{-1} \circ F_x ,
\]
and we trivially have the  identity, cf. \cite[Definition 2.17]{hairer2014theory}.
\begin{equation}\label{gamma-pi}
\Pi_x =  \Pi_y \circ \Gamma_{xy}.
\end{equation}
%
%
%%
%%
%%
%%%As mentioned above, the solutions $\phi$ to \eqref{equation} are constructed in a space of functions $\Phi \colon \R \times \R^3 \to T$, 
%%satisfying a certain continuity condition, the so-called space of modelled distributions of order $\gamma \in  (1 + 2 \epsilon,  \frac32-5 \epsilon)$ , \cite[Definition 3.1]{hairer2014theory}. It turns out 
%%that this function has the form 
%%\begin{equation} \label{def_phi_regstr}
%%\Phi(x) = \oub + v(x) \oneb - \outzb + \nu(x) . \blue{X} - 3v(x) \oudzb,
%%\end{equation}
%%%
%%for some functions  $v$ and $\nu$ (which of course ultimately coincide with our functions $v$ and $\nu$). 
%%
%
 The continuity assumption on $\mathcal{U}$ and $\mathcal{V}$ is  formulated in terms of the translation operators $\Gamma_{xy}$. 
 $\mathcal{U}$ is said to be a modelled distribution of order $\gamma $  if
\begin{equation*}
 \| \mathcal{U}(x) - \Gamma_{xy} \mathcal{U}(y) \|_{\beta} \lesssim d(x,y)^{\gamma - \beta},
\end{equation*}
where $\| \cdot \|_{\beta}$ refers to the component in $T_\beta$. It is easy to check that  for both, $\mathcal{U}$ defined by \eqref{def_phi_regstr}
and $\mathcal{V}$ defined by \eqref{def_V_regstr},  this condition translates precisely into 
the ``modelledness conditions''
\begin{align}
\notag
& | v(y)-v(x)+\outz(y)-\outz(x)- \nu(x) .X(y-x) +  3v(x)(\oudz(y)-\oudz(x)) |  
\lesssim d(x,y)^{\gamma} ,\\
\label{def-norm-modelled-distribution}
&| \nu(y) - \nu(x) | \lesssim d(x,y)^{\gamma -1}, \\
\notag
&|v(y) - v(x)|  \lesssim d(x,y)^{\gamma-1+2\eps},
\end{align}
and this condition for $\gamma = \frac52-3\epsilon$ corresponds exactly to the regularity assumptions on $\mathcal{U}$, $\nu$ and $v$ we work with.

The main feature of the space of modelled distributions is that although expansions like \eqref{def_phi_regstr} are ultimately used as good local 
descriptions of distributions, one can multiply them as if they were of positive regularity, provided one can expand the action of the model to new symbols that are seen as products of the symbols introduced earlier. For equation \eqref{equation} one has to get a bound on $u^3=(v+\ou)^3=v^3+3v^2\ou+3v\oud+\out$. We aim to bound this in terms of:
\begin{itemize}
 \item A high-regularity norm on $v$, namely the $\mathcal{D}^\gamma$ norm of the modelled 
distribution $\mathcal{V}$, which is defined as the smallest possible constant in the inequalities \eqref{def-norm-modelled-distribution};
\item The low regularity $L^\infty$ norm $\| v\|$;
\item The bounds on the various stochastic terms.
\end{itemize} 

The term $v^3$ can immediately be bounded by $\| v \|^3$ and $\out$ is a stochastic term which does not involve $v$. The only terms which require work are
$3v^2\ou$ and $3v\oud$. The distribution $\oud$ has regularity $-1-2\eps$, so a description of $v$ to order $\gamma > 1+2\eps$ is required. Such a description 
is precisely provided by \eqref{def_V_regstr}. One now defines new symbols
\begin{align}\label{symbols-rhs1}
\{\outb,\,   \oudb,  \,  \outdb, \, \ouddb, \,  \oudXb  \},
\end{align} 
associates to them a homogeneity using the rule $| \tau \bar{\tau}| = | \tau| + | \bar{\tau}|$,
and simply defines a new modelled distribution for the local description of $v \oud$ by
\begin{equation} \label{def_Voud_regstr}
\mathcal{V} \oud (x) =  v(x) \oudb + \outdb  - 3v(x) \ouddb- \nu(x) . \oudXb.
\end{equation}
This definition becomes substantial by extending the model $(\Pi_x, \Gamma_{xy})$ to these new symbols. One would like to extend the operator $\Pi_x $ to these products simply 
by defining locally
\[
\Pi_x (\tau \bar{\tau}) (y)=  (\Pi_x \tau (y))  (\Pi_x \bar{\tau} (y)),
\]
but such a definition may not be meaningful when the regularisation is removed.  Fortunately, there is some flexibility at this level.  The main requirements for multiplication to be 
well-behaved are only that \eqref{model-continuity} and \eqref{gamma-pi} remain valid for the new symbols and additionally that one has the identity 
%\comment{mention that the action of the group is determined by its action on $\tau$ and $\bar \tau$}
\[
\Gamma_{xy} (\tau \bar \tau) = (\Gamma_{xy}  \tau)  (\Gamma_{xy}\bar \tau).
\]
It is here that the positive renormalisation, and hence the action of the structure group $G$ becomes more involved than subtracting the value at a base point and the condition $| \tau \bar{\tau}| = | \tau| + | \bar{\tau}|$ becomes strictly stronger than \hol regularity.
For example $\mathbf{\Pi}\outdb$ is a distribution of regularity $-1-2\epsilon$ but its homogeneity is strictly larger, namely $-\frac12-5\epsilon$. The condition \eqref{model-continuity} states that near any base point $x$, $\mathbf{\Pi}\outdb$ is well described by a $\outz(x)\mathbf{\Pi}\oudb$ up to an error of order $-\frac12-5\epsilon$, which is strictly stronger than a bound on the $C^{-1-2\epsilon}$ norm of it.
Our definitions \eqref{def_tree2} and the assumed bounds \eqref{def_outd}, \eqref{def_oudd}, \eqref{def_oudX} correspond exactly to the definitions for $\Pi_x$ and the bound
\eqref{model-continuity} in \cite{hairer2014theory}.
The only difference is that for Hairer defines the trees $\ou$,$\oud$ and $\out$ using the inverse heat operator with some cut-off at large scales and appropriate right hand side. We only assume that they satisfy the heat equation point-wise without imposing any boundary conditions, but we additionally impose some natural regularity bounds, as explained in Remark \ref{The Remark 2}.
 Combining Hairer's multiplication theorem \cite[Theorem 4.7]{hairer2014theory}
and his reconstruction theorem \cite[Theorem 3.10]{hairer2014theory} then yields the estimate
\[
\big| \big\langle \mathcal{R} (\mathcal{V} \oud) -  \Pi_x  (\mathcal{V} \oud), \varphi_x^\lambda \big \rangle \big| \lesssim \lambda^{\gamma-1-2\eps }  \| \mathcal{V} \|_{\mathcal{D}^\gamma} \| \Pi \|, 
\]
where $ \| \Pi \|$ is the smallest possible constant in all of the assumed bounds on the model. This is essentially the statement of our Lemma~\ref{lem_voud} up to a few points:
\begin{itemize}
\item Some of the terms in $\Pi_x (\mathcal{V} \oud) $  are removed from the left hand side of \eqref{eq:lem_voud} and added to the right hand side using the triangle inequality.
\item We prove these estimates ``by hand'' without using the algebraic machinery discussed above and in particular without introducing the group $G$ to organise the various continuity assumptions.
More precisely, Theorem~\ref{Reconstruction}, is a condensed version of \cite[Theorem 3.10]{hairer2014theory} which contains the key analytic estimate, but assumes the output of the algebraic
machinery. In the case, of \eqref{equation} the algebraic manipulations are not too complex and can be done directly quite easily, and that is precisely what we do in the proof of Lemma
\ref{lem_voud}
\item Along the way we keep track of the precise norms needed in each term, rather than compiling them in   $\| \mathcal{V} \|_{\mathcal{D}^\gamma}$  and  $\| \Pi \|$. This added level of detail is important for 
us, especially when determining the exact exponents of each tree appearing in our final estimate \eqref{equation of the theorem}. %\comment{perhaps discuss renormalisation?}
\end{itemize}
The treatment of the term $v^2 \ou$ goes along similar lines. As $\ou$ has better regularity than $\oud$, namely $-\frac12-\eps$, a local description of $v^2$ is only required to order $> \frac12 +\eps$ and this is provided by
\[
\mathcal{V}^2(x) = v^2(x) \oneb -2 v(x) \outzb ,
\]
which in turn prompts us to define
\[
\mathcal{V}^2 \ou (x) = v^2(x) \oub -2 v(x) \outub .
\]
%\comment{good to point out that our bounds only depend linearly on high regularity norm}
Again, our assumption \eqref{def_outu} corresponds exactly to the homogeneity condition \eqref{model-continuity} in \cite{hairer2014theory} and our Lemma~\ref{lem_v2ou} is obtained by combining the 
 multiplication and reconstruction theorem and applying the triangle inequality, this time to remove the term corresponding to  $\Pi_x \mathcal{V}^2 \ou $ from the left hand side completely. The H\"older norm
 $[v+\outz]_{1-2\epsilon}$ which appears on the right hand side of our estimate \eqref{eq:lem_v2ou} corresponds to the norm of the modelled distribution one obtains by removing the terms
 $- 3v(x) \oudzb- \nu(x) . \blue{X}$, which are not necessary here, from the definition of $\mathcal{V}$ in \eqref{def_V_regstr}.

The last ingredient from the theory of regularity structures concerns the heat operator. For us, the gain of regularity for solutions to the heat equation is expressed in Lemma \ref{lemschauder}, and this corresponds to \cite[Theorem 5.12]{hairer2014theory}. Since $\heat^{-1}$ is a $2$-regularising operator ($\beta=2$ in Hairer's theory) and we aim to describe the \lhs of \eqref{eq v} with regularity $\frac32-5\epsilon$, it is enough to describe the \rhs with regularity $-\frac12-5\epsilon$. It is therefore sufficient to work with 
\[
\mathcal{W}(x)=v(x)\oudb+\outb.
\]
Imposing that $\mathcal{W}$ is a modelled distribution of order $-\frac12-5\epsilon$ \cite[Definition 3.1]{hairer2014theory} translates precisely to our three-point continuity condition \eqref{lemschauder2}, and our smallness assumption \eqref{lemschauder1} corresponds to \cite[Equation 5.42]{hairer2014theory}, which in the notation of this section would be
\begin{equation}\label{this equation}
\big{|}\langle \mathcal{R}\mathcal{W}-\Pi_x\mathcal{W},\varphi_x^\lambda\rangle\big{|}\les \lambda^{-\frac12-5\epsilon}.
\end{equation}
Note that the information contained in the modelled distribution $\mathcal{W}$ is not enough to uniquely determine the reconstruction of the right-hand side, which therefore has to be viewed as input for the theorem. The exact statement of \eqref{lemschauder2} appears slightly stronger than \eqref{this equation} because of the $L^{\infty}$ norm on the \lhs and the extra parameter $L$, but in practice the seemingly stronger bound can be obtained easily from the weaker bound using triangle inequality and some lower-order regularity information.

In the framework of regularity structure the operator that encodes the integration of a modelled distribution is described as the sum of three operators. The first operator $\mathcal{I}$ acts point-wise on the modelled distributions by a shift of coefficients. The action on the trees is:
\[
\mathcal{I}\outb=\outzb,\quad\mathcal{I}\oudb=\oudzb.
\]
The continuity of the coefficients for a modelled distribution is transferred accordingly under the action of $\mathcal{I}$. In our setting, this is also automatic and follows from our assumptions, as explained in Remark \ref{The Remark 2}.

The non-trivial part of the integration happens on the levels $\oneb$ and $\blue{X}$ which is encoded in Hairer's theory in the operators $\mathcal{J}$ and $\mathcal{N}$. We have again a direct translation, although we do not need to split the operator.
\begin{align*}
\mathcal{N}\mathcal{W}(x) &= \big( \heat^{-1}(-v^3-3v^2\ou-3(v-v(x))\oud)|_{y=x}\big)  \oneb + \nu(x) \blue{X},\\
 \mathcal{J}(x)\mathcal{W} &= \big( 3v(x)\oudz-\outz  \big)\oneb. 
\end{align*}

The differences between our approach and the one adopted by Hairer is that, in the spirit of \cite{OW} we use a kernel-free approach and we have a special treatment of the boundary on the levels $\oneb$ and $\blue{X}$. We are also more precise in our final bounds in the sense that, as in the definition of $\mathcal{V}\oud$, we keep track of the precise norms needed in each term.

\section{Proof of Theorem \ref{The Theorem}}\label{The Proof}
%The proof relies on two arguments. Once we have made sense of some ill-defined products thanks to the reconstruction theorem,  we control the small scale oscillations via Schauder theory and the large scale behaviour through the maximum principle. Connection between the two is achieved via the convolution of the equation with the kernel introduced in Section~\ref{sec:Main}, which produces a commutator term. The technicality of the proof lies in balancing the contribution of the commutator and the contribution of the irregular noise.

\subsection{Assumption}
We assume that the bound of Theorem \ref{The Theorem} in terms of powers of trees does not hold on a domain $D=P_r$, and use that assumption to prove that then bound in $\frac1R$ holds. Our assumption is stated as such:
\begin{equation}\label{ass1}
\forall\tau\in L,\quad [\tau]_{\gamma_\tau}\leq c\|v\|_{D}^{n_\tau(\frac12-\epsilon)},
\end{equation}
for some constant $c<1$ that we will tune later, according to conditions suggested by equations \eqref{simple schauder} and \eqref{max free}. 
With these assumptions, Lemmas \ref{lem_v2ou} and \ref{lem_voud} can be restated as, for any $x$ with $B(x,T)\in D$,
\begin{align}\label{eq:lem_v2ou_ass}
|(v^2\ou)_T(x)|\les& c([v+\outz]_{1-2\epsilon,B(x,T)}T^{\frac12-3\epsilon}\|v\|_D^{\frac32-\epsilon}+[v]_{\frac12-3\epsilon,B(x,T)}T^{\frac12-7\epsilon}\|v\|_D^{2-4\epsilon})\nonumber\\
&+c(T^{-\frac12-\epsilon}\|v\|_D^{\frac52-\epsilon}+ T^{-4\epsilon}\|v\|_D^{3-4\epsilon}+T^{\frac12-7\epsilon}\|v\|_D^{\frac72-7\epsilon}).
\end{align}
and:
\begin{align}\label{eq:lem_voud_ass}
|((v&-v(x))\oud)_T(x)+3C_2(v_T+\ou_T)(x)| \nonumber\\
\les&cT^{\frac12-7\epsilon}\Big([v]_{\frac12-3\epsilon,B(x,T)}\|v\|_D^{2-4\epsilon}
+([U]_{\frac32-5\epsilon,B(x,T)} +[\nu]_{\frac12-5\epsilon,B(x,T)})\|v\|_D^{1-2\epsilon}\Big)\nonumber\\
&+c(\|v\|_{D}^{\frac52-5\epsilon}T^{-\frac12-5\epsilon}+T^{-4\epsilon}\|v\|_D^{3-4\epsilon})+\|\nu\|_{B(x,T)}c\|v\|_D^{1-2\epsilon}T^{-2\epsilon}.
\end{align}

\subsection{Applying Lemma \ref{lemschauder}}

For any domain $D$, for $x,L,T$ and $d$ such that $B(x,L+T)\subset D_d$, we prove the following bound, which is the first condition to apply Lemma~\ref{lemschauder}.
\begin{align}\label{bound heat U}
\|\heat &U(x,\cdot)_T\|_{B(x,L)}\les\|v\|_{D}^3+ L^{\frac12-3\epsilon}T^{-1-2\epsilon}c[v]_{\frac12-3\epsilon,D_d,d}\|v\|_{D}^{1-2\epsilon}\nonumber\\
+c&([v+\outz]_{1-2\epsilon,B(x,T)}T^{\frac12-3\epsilon}\|v\|_D^{\frac32-\epsilon}+[v]_{\frac12-3\epsilon,B(x,T)}T^{\frac12-7\epsilon}\|v\|_D^{2-4\epsilon})\nonumber\\
+c&(T^{-\frac12-\epsilon}\|v\|_D^{\frac52-\epsilon}+ T^{-4\epsilon}\|v\|_D^{3-4\epsilon}+T^{\frac12-7\epsilon}\|v\|_D^{\frac72-7\epsilon})\\
+c&T^{\frac12-7\epsilon}\Big([v]_{\frac12-3\epsilon,B(x,T)}\|v\|_D^{2-4\epsilon}
+([U]_{\frac32-5\epsilon,B(x,T)} +[\nu]_{\frac12-5\epsilon,B(x,T)})\|v\|_D^{1-2\epsilon}\Big)\nonumber\\
+c&(\|v\|_{D}^{\frac52-5\epsilon}T^{-\frac12-5\epsilon}+T^{-4\epsilon}\|v\|_D^{3-4\epsilon})+\|\nu\|_{B(x,T)}c\|v\|_D^{1-2\epsilon}T^{-2\epsilon},\nonumber
\end{align}
where $\nu$ is the optimal function in the definition of $[U]_{\frac32-5\epsilon,D_d}$.

Let $x$ be an arbitrary point in $D_{d+L}$ and $y$  a point in $B(x,L)\subset D_d$. We have
\begin{align*}
\heat &U(x,\cdot)_T(y)=\int\Psi_T(z-y)\heat U(x,z)dz\\
=&-(v^3)_T(y)-3(v^2\ou)_T(y)-3((v-v(x))\oud)_T(y)-9C_2(v_T(y)+\ou_T(y))\\
=&-(v^3)_T(y)-3(v^2\ou)_T(y)-3(v(y)-v(x))\oud_T(y)\\&-3((v-v(y))\oud)_T(y)-9C_2(v_T(y)+\ou_T(y)).
\end{align*}
We bound the some terms of this expression by the previous bounds \eqref{eq:lem_v2ou_ass} and \eqref{eq:lem_voud_ass} and the remaining ones as follows:
\begin{align}\label{different terms}
|(v^3)_T(y)|&\leq\|v\|_{B(y,T)}^3\leq \|v\|_{D_d}^3 ,\\
|(v(y)-v(x))\oud_T(y)|&\leq d(x,y)^{\frac12-3\epsilon}[v]_{\frac12-3\epsilon,D_d,d}T^{-1-2\epsilon}[\oud]_{-1-2\epsilon}\nonumber\\
&\leq L^{\frac12-3\epsilon}T^{-1-2\epsilon}c[v]_{\frac12-3\epsilon,D_d,d}\|v\|_{D}^{1-2\epsilon}\label{different terms bis}
\end{align}
This proves \eqref{bound heat U}.

The three-point continuity on $U$ holds as follows. for any $x\in D_d$, for any $y\in B(x,\frac{d}4)$, for any $z\in B(y,\frac{d}4)$
\begin{align}\label{3pc}
|U(x,y)-U(x,z)-U(z,y)|=& 3|v(x)-v(z)||\oudz(y)-\oudz(z)|\nonumber\\
\leq &3[v]_{\frac12-3\epsilon,D_\frac{d}2,\frac{d}2}[\oudz]_{1-2\epsilon}d(x,z)^{\frac12-3\epsilon}d(y,z)^{1-2\epsilon}\\
\leq &3c[v]_{\frac12-3\epsilon,D_\frac{d}2,\frac{d}2}\|v\|_{D}^{1-2\epsilon}d(x,z)^{\frac12-3\epsilon}d(y,z)^{1-2\epsilon}\nonumber	
\end{align}

Lemma \ref{lemschauder} applies to $U$ with $\kappa=\frac32-5\epsilon$. Note that in the bound \ref{bound heat U} we see powers of $T$ higher than $T^{-\frac12-5\epsilon}$ but we use the fact that $T\les d$ to make up for that.
 After a few simplifications we get

\begin{align}\label{apply schauder}
\sup_{d\leq d_0} &d^{\frac32-5\epsilon}[U]_{\frac32-5\epsilon,D_d}\les\sup_{d\leq d_0} \Big(d^2\|v\|_{D}^3+ d^{\frac32-5\epsilon}c[v]_{\frac12-3\epsilon,D_d,d}\|v\|_{D}^{1-2\epsilon}\nonumber\\
&+c([v+\outz]_{1-2\epsilon,B(x,T)}d^{\frac52-3\epsilon}\|v\|_D^{\frac32-\epsilon}+d^{\frac32-\epsilon}\|v\|_D^{\frac52-\epsilon}+d^{\frac52-7\epsilon}\|v\|_D^{\frac72-7\epsilon})\nonumber\\
&+cd^{\frac52-7\epsilon}\Big([v]_{\frac12-3\epsilon,B(x,T)}\|v\|_D^{2-4\epsilon}
+([U]_{\frac32-5\epsilon,B(x,T)} +[\nu]_{\frac12-5\epsilon,B(x,T)})\|v\|_D^{1-2\epsilon}\Big)\\
&+c(\|v\|_{D}^{\frac52-5\epsilon}d^{\frac32-5\epsilon}+d^{2-4\epsilon}\|v\|_D^{3-4\epsilon})+\|\nu\|_{B(x,T)}c\|v\|_D^{1-2\epsilon}d^{2-2\epsilon}+\|U\|_{D_d,d}\Big).\nonumber	
\end{align}

\subsection{Simplifications}

Our goal in this section is to produce bounds on the semi-norms $[v]_{\frac12-3\epsilon,D_d,d},[v+\outz]_{1-2\epsilon,D_d,d}$ and $[U]_{\frac32-5\epsilon,D_d}$ that depend only on $\|v\|_{D_d}$, in particular independent of each other. 
We introduce the following elementary bounds, which can be deduced from triangle inequalities and application of the Assumption \ref{ass1}.
\begin{align}
[v]_{\frac12-3\epsilon,D_d,d}&\leq d^{1-2\epsilon}[U]_{\frac32-5\epsilon,D_d}+[\outz]_{\frac12-3\epsilon}+3d^{\frac12+\epsilon}\|v\|_{D_d}[\oudz]_{1-2\epsilon}+d^{\frac12+3\epsilon}\|\nu\|_{D_d}\nonumber\\
& \leq d^{1-2\epsilon}[U]_{\frac32-5\epsilon,D_d}+c\|v\|_{D}^{\frac32-3\epsilon}+3cd^{\frac12+\epsilon}\|v\|_{D_d}^{2-2\epsilon}+d^{\frac12+3\epsilon}\|\nu\|_{D_d} \label{ele est v}
\end{align}
and
\begin{align}
[v+\outz]_{1-2\epsilon,D_d,d}&\leq d^{\frac12-3\epsilon}[U]_{\frac32-5\epsilon,D_d}+3\|v\|_{D_d}[\oudz]_{1-2\epsilon}+d^{2\epsilon}\|\nu\|_{D_d}\nonumber\\
&\leq d^{\frac12-3\epsilon}[U]_{\frac32-5\epsilon,D_d}+3c\|v\|_{D_d}^{2-2\epsilon}+d^{2\epsilon}\|\nu\|_{D_d}\label{ele est v outz},
\end{align}
and we recall that from Corollary \ref{corschauder} we have the two bounds, assuming $d\in (0,r_d]$,
\begin{align}
\|\nu\|_{D_d}\les&[U]_{\frac32-5\epsilon,D_d}d^{\frac12-5\epsilon}+\|U\|_{D_d,d}d^{-1}\label{ele est nu}
\end{align}
and
\begin{align}
[\nu]_{\frac12-5\epsilon,D_d}\les& [U]_{\frac32-5\epsilon,D_d}+[v]_{\frac12-3\epsilon,D_d,d}[\oudz]_{1-2\epsilon}+d^{-\frac32+5\epsilon}\|U\|_{D_d,d}\nonumber\\
\les& [U]_{\frac32-5\epsilon,D_d}+c[v]_{\frac12-3\epsilon,D_d,d}\|v\|_D^{1-2\epsilon}+d^{-\frac32+5\epsilon}\|U\|_{D_d,d}\label{ele est nu hol}.
\end{align}
We will also be using the bound 
\begin{align}\label{ele est U}
\|U\|_{D_d,d}\leq& 2\|v\|_{D}+[\outz]_{\frac12-3\epsilon}d^{\frac12-3\epsilon}+3\|v\|_{D}[\oudz]_{1-2\epsilon}d^{1-2\epsilon}\nonumber\\
\leq& 2\|v\|_{D}+cd^{\frac12-3\epsilon}\|v\|_{D}^{\frac32-3\epsilon}+3cd^{1-2\epsilon}\|v\|_{D}^{2-2\epsilon}.
\end{align}
By combining the bounds above to get bounds only in terms of $\|v\|_{D}$ and  $[U]_{\frac32-5\epsilon,D_d}$ we get the following bounds (in order of logical deduction):
\begin{align}
&\|\nu\|_{D_d}\les d^{\frac12-5\epsilon}[U]_{\frac32-5\epsilon,D_d}+d^{-1}\|v\|_{D_d}+cd^{-\frac12-3\epsilon}\|v\|_{D_d}^{\frac32-3\epsilon}+cd^{-2\epsilon}\|v\|_{D_d}^{2-2\epsilon},\label{comb1}\\
&[v+\outz]_{1-2\epsilon,D_d,d}\les d^{\frac12-3\epsilon}[U]_{\frac32-5\epsilon,D_d}+c\|v\|_{D_d}^{2-2\epsilon}
+d^{-1+2\epsilon}\|v\|_{D_d}\nonumber\\&\qquad\qquad\qquad\qquad+cd^{-\frac12-\epsilon}\|v\|_{D_d}^{\frac32-3\epsilon},\\
&[v]_{\frac12-3\epsilon,D_d,d}\les d^{1-2\epsilon}[U]_{\frac32-5\epsilon,D_d}+d^{-\frac12+3\epsilon}\|v\|_{D_d}+c\|v\|_{D_d}^{\frac32-3\epsilon}\label{comb3}\\&\qquad\qquad\qquad+cd^{\frac12+\epsilon}\|v\|_{D_d}^{2-2\epsilon},\nonumber\\
&[\nu]_{\frac12-5\epsilon,D_d}\les [U]_{\frac32-5\epsilon,D_d}(1+cd^{1-2\epsilon}\|v\|_D^{1-2\epsilon})
+c^2(\|v\|_{D_d}^{\frac52-5\epsilon}+d^{\frac12+\epsilon}\|v\|_{D_d}^{3-4\epsilon})\nonumber\\
&\qquad\qquad\qquad+d^{-\frac32+5\epsilon}\|v\|_{D_d}+cd^{-1+2\epsilon}\|v\|_{D_d}^{\frac32-3\epsilon}+cd^{-\frac12+3\epsilon}\|v\|_{D_d}^{2-2\epsilon}.\label{comb4}
\end{align}

We inject those in the \rhs of \eqref{apply schauder} and we bound positive powers of $d$ by powers of $d_0$. A few computations allow to reduce the result to the following expression:
\begin{align}\label{simple schauder}
\sup_{d\leq d_0}& d^{\frac32-5\epsilon}[U]_{\frac32-5\epsilon,D_d}
%\les\sup_{d\leq d_0} \Big(d^2\|v\|_{D}^3\\
%+ cd^{1-2\epsilon}\|v\|_{D}^{2-2\epsilon}+(c+c^2)(d^{\frac32-5\epsilon}\|v\|_{D}^{\frac52-5\epsilon}+d^{2-4\epsilon}\|v\|_{D}^{3-4\epsilon})\nonumber\\
%+c^2(d^{\frac52-3\epsilon}\|v\|_D^{\frac32-\epsilon}+(1+c)d^{3-6\epsilon}\|v\|_D^{2-4\epsilon})( \|v\|_{D_d}^{2-2\epsilon}+d^{-\frac12-\epsilon}\|v\|_{D_d}^{\frac32-3\epsilon})\\
%+c(d^{\frac32-\epsilon}\|v\|_D^{\frac52-\epsilon}+ (1+c)d^{2-4\epsilon}\|v\|_D^{3-4\epsilon}+(1+c)d^{\frac52-7\epsilon}\|v\|_D^{\frac72-7\epsilon})\\
%+[U]_{\frac32-5\epsilon,B(x,T)} (cd^{\frac52-7\epsilon}\|v\|_D^{1-2\epsilon}+c(d^{3-6\epsilon}\|v\|_D^{\frac32-\epsilon}+(1+c)d^{\frac72-9\epsilon}\|v\|_D^{2-4\epsilon}))\nonumber\\
%+c^3d^{\frac52-7\epsilon}(\|v\|_{D_d}^{\frac72-7\epsilon}+cd^{3-6\epsilon}\|v\|_{D_d}^{\frac{11}2-11\epsilon})\nonumber\\
%+\|v\|_{D}+cd^{\frac12-3\epsilon}\|v\|_{D}^{\frac32-3\epsilon}\Big).\nonumber	\\
\les c\sum_{h}d_0^h\|v\|_{D}^{h+1}
+c\sup_{d\leq d_0}d^{\frac32-5\epsilon}[U]_{\frac32-5\epsilon,D_d}\sum_{l}d_0^l\|v\|_{D}^l,
\end{align}
where the index of the sum $h$ is taken in a finite subset of $[0,\frac92-11\epsilon]$ and the index $l$ in $\{1-2\epsilon,\frac23-\epsilon,2-4\epsilon\}$. 
If we set 
\begin{equation}\label{def d0}
d_0=\|v\|_D^{-1},
\end{equation} 
 we see that we can get rid of $\sup_{d\leq d_0}d^{\frac32-5\epsilon}[U]_{\frac32-5\epsilon,D_d}$ in the \rhs under a first smallness condition on $c$ depending on the constant implicit in $\les$. If this condition is satisfied, we have:
\begin{align}\label{simplest schauder}
\sup_{d\leq d_0} d^{\frac32-5\epsilon}[U]_{\frac32-5\epsilon,D_d}\les c\|v\|_D.
\end{align}
In this equation and in the following, $\les$ does not depend on this first condition on $c$.
Applying this to Equations \eqref{comb1} to \eqref{comb4} gives 
\begin{align}
\sup_{d\leq d_0}d\|\nu\|_{D_d}\les c\|v\|_D,\label{simplest schauder nu}\\
\sup_{d\leq d_0}d^{1-2\epsilon}[v+\outz]_{1-2\epsilon,D_d,d}\les c\|v\|_D,\label{simplest schauder v+outu}\\
\sup_{d\leq d_0}d^{\frac12-3\epsilon}[v]_{\frac12-3\epsilon,D_d,d}\les c\|v\|_D,\label{simplest schauder v}\\
\sup_{d\leq d_0}d^{\frac32-5\epsilon}[\nu]_{\frac12-5\epsilon,D_d}\les c\|v\|_D.\label{simplest schauder nu hol}
\end{align}

Applying estimates \eqref{simplest schauder} and \eqref{simplest schauder nu}-\eqref{simplest schauder nu hol} to \eqref{eq:lem_v2ou_ass} and \eqref{eq:lem_voud_ass}, we have, for $d\leq d_0$ and $B(x,T)\in D_d$,
\begin{align}\label{eq:lem_v2ou_ass_sh}
|(v^2\ou)_T(x)|\les& c^2(d^{-1+2\epsilon}T^{\frac12-3\epsilon}\|v\|_D^{\frac52-\epsilon}+d^{-\frac12+3\epsilon}T^{\frac12-7\epsilon}\|v\|_D^{3-4\epsilon})\\
&+c(T^{-\frac12-\epsilon}\|v\|_D^{\frac52-\epsilon}+ T^{-4\epsilon}\|v\|_D^{3-4\epsilon}+T^{\frac12-7\epsilon}\|v\|_D^{\frac72-7\epsilon})\nonumber.
\end{align}
and
\begin{align}\label{eq:lem_voud_ass_sh}
|(&v\oud)_T(x)+3C_2(v_T+\ou_T)(x)|\les c\|v\|_D^{2-2\epsilon}T^{-1-2\epsilon}\nonumber\\
&+c^2T^{\frac12-7\epsilon}\Big(d^{-\frac12+3\epsilon}\|v\|_D^{3-4\epsilon}
+d^{-\frac32+5\epsilon}\|v\|_D^{2-2\epsilon}\Big)\nonumber\\
&+c(\|v\|_{D}^{\frac52-5\epsilon}T^{-\frac12-5\epsilon}+T^{-4\epsilon}\|v\|_D^{3-4\epsilon})+c^2d^{-1}\|v\|_D^{2-2\epsilon}T^{-2\epsilon}.
\end{align}
In this last estimate, we have used triangle inequality to get $\|v\oud_T\|_{D_{d+T}}$ out of the left hand side, and then used the assumption \ref{ass1} to bound it.
%We will apply this in the following section with $D_d=P_{R'-T}$, $d=(\lambda -1) T\sim\lambda  T=d_0$ for some $\lambda >2$ to be defined, hence $D=P_{R'-\lambda  T}$.\comment{need to move this in a sensible place}

\subsection{Application of Lemma \ref{main thm poly 1}}
We now go back to the original equation, and start to study large scale.
We convolve the equation \eqref{eq v} with $\Psi_T$:
\begin{equation}\label{phi42_T}
(\partial_t-\Delta)(v)_T=-(v_T)^3-3(v^2\ou)_T-3(v\oud)_T-9C_2(v_T+\ou_T)-(\out)_T+((v_T)^3-(v^3)_T).
\end{equation}
Lemma \ref{main thm poly 1} implies that for all $r>0$ and $0<R'<R$ such that $r+R<\frac12$, we have

\begin{align}\label{apply max}
\|(v)_T\|_{P_{r+R}}\les\max&\Big{\{}\frac1{R-R'},\|(v_T)^3-(v^3)_T\|_{P_{r+R'}}^\frac13,\|(v^2\ou)_T\|_{P_{r+R'}}^\frac13,\nonumber\\
&\|(v\oud)_T+3C_2(v_T+\ou_T)\|_{P_{r+R'}}^\frac13,\|(\out)_T\|_{P_{r+R'}}^\frac13\Big{\}}.
\end{align}
The goal is now to balance the commutator and the renormalized powers of the noise term by choosing the parameter $T$ appropriately. We first mention that applying \eqref{mollified regularity 1} and then Assumption \ref{ass1} gives the bound:
\begin{align}\label{apply max bis}
\|v\|_{P_{r+R}}\les\max&\Big{\{}\frac1{R-R'},T^{\frac12-3\epsilon}[v]_{\frac12-3\epsilon,P_{r+R-T},2T},\|(v_T)^3-(v^3)_T\|_{P_{r+R'}}^\frac13,\nonumber\\
&\|(v^2\ou)_T\|_{P_{r+R'}}^\frac13,\|(v\oud)_T+3C_2(v_T+\ou_T)\|_{P_{r+R'}}^\frac13,\|(\out)_T\|_{P_{r+R'}}^\frac13\Big{\}}.
\end{align}
 We need estimates on the commutator $(v_T)^3-(v^3)_T$. This is easily obtained as $v$ is $C^{\frac12-3\epsilon}$, using the moment bounds \eqref{moment of psi} and \eqref{mollified regularity 1}. For any $z\in P_{r+R'}$,
\begin{align*} 
((v&_T)^3-(v^3)_T)(z)=\int\Psi_T(z-\bar{z})\left(v_T(z)^3-v(\bar{z})^3\right)d\bar{z}\\
=&\int\Psi_T(z-\bar{z})\int_0^1\left(v_T(z)-v(\bar{z})\right)3\left(\lambda v_T(z)+(1-\lambda )v(\bar{z})\right)^{2}d\lambda  d\bar{z}\\
\leqslant& 3\|v\|_{B(z,T)}^{2}\int\Psi_T(z-\bar{z})\left(v_T(z)-v(z)+v(z)-v(\bar{z})\right)d\bar{z}\\
\leqslant& 3\|v\|_{B(z,T)}^{2}\int\Psi_T(z-\bar{z})\left(T^{\frac12-3\epsilon}[v]_{\frac12-3\epsilon,B(z,T)}+[v]_{\frac12-3\epsilon,B(z,T)}d(z,\bar{z})^{\frac12-3\epsilon}\right)d\bar{z}\\
\leqslant& 6\|v\|_{B(z,T)}^{2}T^{\frac12-3\epsilon}[v]_{\frac12-3\epsilon,B(z,T)}.
\end{align*}
Since this is true for all $z\in B(z,T)$, 
\begin{equation}\label{commutator estimate 1}
\|(v_T)^3-(v^3)_T\|_{P_{r+R'}}\les \|v\|_{P_{r+R'- T}}^2T^{\frac12-3\epsilon}[v]_{\frac12-3\epsilon,P_{r+R'- T},2T}.
\end{equation}
In conclusion of this step,
\begin{align}\label{apply max ter}
\|v\|_{P_{r+R}}&\les\nonumber\\
\max\Big{\{}&\frac1{R-R'},T^{\frac12-3\epsilon}[v]_{\frac12-3\epsilon,P_{r+R-T},2T},T^{\frac16-\epsilon}\|v\|_{P_{r+R'- T}}^\frac23[v]_{\frac12-3\epsilon,P_{r+R'- T},2T}^\frac13,\nonumber\\
&\|(v^2\ou)_T\|_{P_{r+R'}}^\frac13,\|(v\oud)_T+3C_2(v_T+\ou_T)\|_{P_{r+R'}}^\frac13,\|(\out)_T\|_{P_{r+R'}}^\frac13\Big{\}}.
\end{align}

\subsection{Choice of scale}
We now apply the assumption \eqref{ass1} with $\tau=\out$ and the results from the previous steps, \eqref{simplest schauder v},\eqref{eq:lem_v2ou_ass_sh} and \eqref{eq:lem_voud_ass_sh} to equation \eqref{apply max ter}. In \eqref{apply max ter} we choose 
\[
R'=d_0\text{ and }T=\frac{d_0}k,
\]
for some $k>2$ to be specified. Recall that $d_0=\frac1{\|v\|_{P_r}}$, as set in \eqref{def d0}. In the \lhs of \eqref{simplest schauder v},\eqref{eq:lem_v2ou_ass_sh} and \eqref{eq:lem_voud_ass_sh} we make the particular choice 
\[
d=d_0\frac{k-1}k.
\]
Since $k>2$ we have $d\sim d_0$, and we also have $T+d=d_0$ so $\|v\|_{P_{r+R'- T}}=\|v\|_{D_d}$.
 Equation \eqref{apply max ter} simplifies to
\begin{align}\label{max free}
\|v\|_{P_{r+R}}\leqslant C\max&\Big{\{}\frac1{R-R'},\|v\|_{P_r}k^{-\frac12+3\epsilon},\|v\|_{P_r}k^{-\frac16+\epsilon},\nonumber\\
&\|v\|_{P_r}\Big[c^2\Big(k^{-\frac12+3\epsilon}+k^{-1+6\epsilon}+k^{-\frac12+7\epsilon}+k^{-\frac12+7\epsilon}+k^{2\epsilon}\Big)\\
&+c\Big(k^{-\frac12+\epsilon}+ k^{4\epsilon}+k^{-\frac12+7\epsilon}+k^{1+\epsilon}+k^{\frac12+5\epsilon}+k^{\frac32+3\epsilon}\Big)\Big]^\frac13\Big{\}},\nonumber
\end{align}
for some constant $C>1$.
We see that we can choose $k$ large and then impose another smallness condition on $c$ to get, 
\begin{equation}\label{max freest}
\|v\|_{P_{R+r}}\leqslant \max\Big{\{}\frac{C}{R-R'},\frac12\|v\|_{P_r}\Big{\}}.
\end{equation}

\subsection{ Iterating the result} 
If we have $R\geq 2R'$, then we can rewrite equation \eqref{max freest} for $r=0$ as
\begin{equation}\label{max freeest}
\|v\|_{P_R}\leqslant \max\Big{\{}\frac{2C}{R},\frac{\|v\|_{P}}2 \Big{\}}.
\end{equation}
The first argument of the maximum \eqref{max freeest} is equal to the second one for 
 \[R= R_1:=\frac{4C}{\|v\|_{P}}.\]
This is not in contradiction with $R\geq 2R'=\frac{2}{\|v\|_{P}}$ as $C>1$.
We now define a finite set $0=R_0< \ldots  <R_N =\frac12$ by setting   
\[R_{n+1}-R_n=4C\|v\|_{P_{R_{n}}}^{-1},\]
as long as the times $R_{n+1}$ defined this way stay strictly less than $\frac12$. We terminate the sequence once this algorithm would produce 
a $R_{n+1} \geq \frac12$ in which case we set $R_{n+1} = R_N =\frac12$ or once the Assumption \eqref{ass1} does not hold for $D=R_n$. Note that $4C\|v\|_{P_{R_{n}}}^{-1}$  is increasing in $n$ so the sequence necessarily 
terminates after finitely many steps.
Equation \eqref{max freest} applied with $r=R_{n-1}$ for $n=1...N$ then give the bounds for smaller and smaller parabolic boxes.
\begin{equation}
\|v\|_{P_{R+R_{n-1}}}\leqslant \max\Big{\{}\frac{2C}R,\frac12\|v\|_{P_{R_{n-1}}}\Big{\}},
\end{equation}
hence for $R=R_n-R_{n-1}$, 
\begin{equation}\label{equation with Rn}
\|v\|_{P_{R_n}}\leqslant \frac12\|v\|_{P_{R_{n-1}}}.
\end{equation}
We now show that the bound \eqref{equation of the theorem} in Theorem \ref{The Theorem} holds for all $R=R_n, n\in \{0,...,N\}.$ If Assumption \eqref{ass1} does not hold for $D=R_N$ it is immediate, in the other case for $k\leqslant n$, 
  $\|v\|_{P_{R_n}}\leqslant \|v\|_{P_{R_k}}2^{k-n}$ and hence
\begin{align}
R_n=\sum_{k=0}^{n-1}R_{k+1}-R_k=\sum_{k=0}^{n-1}4C\|v\|_{P_{R_{k}}}^{-1}\leqslant4C\|v\|_{P_{R_n}}^{-1}\sum_{k=0}^{n-1}2^{k-n}
\les&\label{e:ppp} \|v\|_{P_{R_n}}^{-1}.
\end{align}
For the end point $R_N$ we have  either $R_{N-1} \geqslant \frac14$ or $R_{N}-R_{N-1}\geqslant \frac14$. In the first case 
we invoke \eqref{e:ppp} for $n = N-1$ and in the second case the definition of $R_{n+1}-R_n$, in both cases yielding a bound on $\|v\|_{P_{R_{N-1}}}$.
Finally for values $R\in(R_n,R_{n+1})$, we use the definition of $R_{n+1}-R_n$: 
\begin{align*}
R&\leqslant R_{n+1}= R_{n+1}-R_n+R_n
\les \|v\|_{P_{R_n}}^{-1}+R_n\leq \|v\|_{P_R}^{-1}.
\end{align*}
This concludes the proof of the theorem.

\section{Proof of the lemmas}\label{proofproof}
\subsection{Proof of Lemma \ref{lemschauder}}\label{proof schauder}
This proof is inspired from the one in \cite[Proposition 2]{otto2018parabolic}, our contribution being the introduction of blow-up at the boundaries of the domain instead of an assumption of periodicity.

\paragraph*{\textsc{Step} 1.}

We claim that for all base points $x$ and scales $T,R$ and $L$
with $R\leq \frac{L}2$ and such that $B(x,L)\subset D$, it holds:
\begin{align}\label{Schauder Step 1}
&\inf_l\|U_T(x,.)-l\|_{B(x,R)} \\
\notag
&\les\frac{R^2}{L^2}\inf_l\|U_T(x,.)-l\|_{B(x,L)}+L^2M^{(1)}_{\{x\},L}\sum_{\beta\in A}T^{\beta-2} L^{\kappa-\beta},
\end{align}
where the infimum runs over all affine functions $l(y)=C.X(y-x)+c$.
To prove this, we define a decomposition $U_T(x,.)=u_>+u_<$ where $u_>$ is the solution to 
\[
\heat u_>=\textbf{1}_{B(x,L)}\heat U_T(x,.).
\]
with Dirichlet boundary conditions. By standard estimates for the heat equation \cite[Cor.8.1.5]{krylov1996lectures},
\begin{equation}\label{sh2int}
\|u_>\|_{B(x,L)}\les L^2\|\heat U_T(x,\cdot)\|_{B(x,L)}\overset{\eqref{lemschauder1}}{\leqslant}L^2M^{(1)}_{\{x\},L}\sum_{\beta\in A}T^{\beta-2} L^{\kappa-\beta}.
\end{equation}
As $(\partial_t-\Delta)u_<=0$ on $B(x,L)$ for $\partial\in\{\partial_t,\partial_{i} \partial_{j}\}$ a differential operator of order $1$ in time or $2$ in space,
\[
\|\partial u_<\|_{B(x,R)}\leqslant L^{-2}\|u_<-l_>\|_{B(x,L)},
\]
for any affine function $l_>$, where we used $R\leqslant \frac{L}2$, and the fact that the differential operators used cancel the spatial linear functional. Next we define a concrete affine function $l_<$ via $l_<(y) := u_<(x) + \grad u_<(x).X(y - x)$ and observe, using Taylor's formula,
\begin{align*}
\|u_<-l_<\|_{B(x,R)}\leq &R^2\|Du_<\|_{B(x,R)}\\
\leq&\frac{R^2}{L^2}\|u_<-l_>\|_{B(x,L)}\\
\leqslant&\frac{R^2}{L^2}\|U_T(x,.)-l_>\|_{B(x,L)}+\|u_>\|_{B(x,L)}.
\end{align*}
Using the triangle inequality once more and \eqref{sh2int} gives:
\begin{align*}
\|U_T(x,.)-l_>\|_{B(x,R)}\leqslant&\|u_<-l_<\|_{B(x,R)}+\|u_>\|_{B(x,R)}\\
\les&\frac{R^2}{L^2}\|U_T(x,.)-l_>\|_{B(x,L)}+\|u_>\|_{B(x,L)}\\
\overset{\eqref{sh2int}}{\leq}&\frac{R^2}{L^2}\|U_T(x,.)-l_>\|_{B(x,L)}+L^2M^{(1)}_{\{x\},L}\sum_{\beta\in A}T^{\beta-2} L^{\kappa-\beta},
\end{align*}
which implies \eqref{Schauder Step 2}

\paragraph*{\textsc{Step} 2.} We claim that for all base points $x$ and scales $T,L$, it holds:
\begin{equation}\label{Schauder Step 2}
\|U_T(x,.)-U(x,.)\|_{B(x,R)}\leq M^{(2)}_{\{x\},R,T}\sum_{\beta\in A} R^\beta T^{\kappa-\beta}+ T^\kappa[U]_{\kappa,B(x,R),T}.
\end{equation}
Indeed, since $\Psi$ is symmetric, it integrates to $0$ against linear functions hence for any $y\in B(x,R)$, we have
\begin{align*}
|U_T(x,y)&-U(x,y)|=\Big|\int \Psi_T(y-z)(U(x,z)-U(x,y))dz \Big|\\
=& \inf_{\nu(y)}\Big|\int \Psi_T(y-z)(U(x,z)-U(x,y)-U(y,z))dz
\\&+\int \Psi_T(y-z)(U(y,z)-\nu(y).X(z-y))dz \Big|\\
\leqslant& M^{(2)}_{\{x\},R,T}\sum_{\beta\in A}d(y,x)^\beta\int \Psi_T(y-z) d(z, y)^{\kappa-\beta}dz\\
+&\Big(\sup_{y\in B(x,R)}\inf_{\nu(y)}\sup_{z\in B(y,T)}d(y,z)^{-\kappa}|U(y,z)-\nu(y).X(z-y)|\Big)\\&\times\int \Psi_T(y-z)d(z,y)^\kappa dz.
\end{align*}

\paragraph*{\textsc{Step} 3.}
We prove \begin{align}\label{app_st_2}
\sup_{R\leq \frac{\epsilon d}2}R^{-\kappa}&\inf_l\|U(x,.)-l\|_{B(x,R)}\nonumber\\
\les \sum_{\beta\in A}&\Big(M^{(1)}_{\{x\},\frac{d}2}\epsilon^{-4+2\beta-\kappa}+M^{(2)}_{\{x\},\frac{\epsilon d}2,\frac{\epsilon^2 d}2} \epsilon^{\kappa-\beta}+M^{(2)}_{\{x\},\frac{ d}2,\frac{\epsilon^2 d}2} \epsilon^{2(\kappa-\beta)}\Big)\nonumber\\
&+\epsilon^{2-2\kappa}\frac{ d^{-\kappa}}{2^{-\kappa}}\|U(x,.)\|_{B(x,\frac{d}2(1+\epsilon^2))}+ (\epsilon^\kappa+\epsilon^{2+\kappa})[U]_{\kappa, B(x,\frac{\epsilon d}2),\frac{\epsilon^2 d}2}.
\end{align}
Multiplying Equation \eqref{Schauder Step 1} by $R^{-\kappa}$ and fixing the length ratios $R=\epsilon L=\epsilon^{-1}T$ for some $\epsilon \leq \frac12$ to be fixed below, we get for any point $x\in D_d$ and length $L\leq \frac{ d}{2}$,

\begin{align*}
&R^{-\kappa}\inf_l\|U_{T}(x,.)-l\|_{B(x,R)}\\&\les \epsilon^{2-\kappa}L^{-\kappa}\inf_l\|U_{T}(x,.)-l\|_{B(x,L)}+\sum_{\beta\in A}M^{(1)}_{D_d,L}\epsilon^{-4+2\beta-\kappa}.
\end{align*}
Taking the supremum over $L\leq \frac{ d}{2}$ while keeping the ratios $R=\epsilon L=\epsilon^{-1}T$ fixed we get
\begin{align*}
&\sup_{R\leq \frac{\epsilon d}2}R^{-\kappa}\inf_l\|U_{T}(x,.)-l\|_{B(x,R)}\\
&\les \epsilon^{2-\kappa}\sup_{L\leq \frac{d}2}L^{-\kappa}\inf_l\|U_{T}(x,.)-l\|_{B(x,L)}+\sup_{L\leq \frac{d}2}\sum_{\beta\in A}M^{(1)}_{D_d,L}\epsilon^{-4+2\beta-\kappa}\\
&\leq  \epsilon^{2-\kappa}\sup_{L\leq \frac{\epsilon d}2}L^{-\kappa}\inf_l\|U_{T}(x,.)-l\|_{B(x,L)}+\sup_{L\leq \frac{d}2}\sum_{\beta\in A}M^{(1)}_{D_d,L}\epsilon^{-4+2\beta-\kappa}\\
&+\epsilon^{2-\kappa}\sup_{\frac{\epsilon d}2\leq L\leq \frac{d}2}L^{-\kappa}\inf_l\|U_{T}(x,.)-l\|_{B(x,L)}.
\end{align*}
The last term is bounded by \[
\epsilon^{2-\kappa}\Big(\frac{\epsilon d}2\Big)^{-\kappa}\|U_{T}(x,.)\|_{B(x,\frac{d}2)}\leq \epsilon^{2-2\kappa}\frac{ d^{-\kappa}}{2^{-\kappa}}\|U(x,.)\|_{B(x,\frac{d}2(1+\epsilon^2))}.\]
 Hence we have
\begin{align*}
\sup_{R\leq \frac{\epsilon d}2}R^{-\kappa}&\inf_l\|U_{T}(x,.)-l\|_{B(x,R)}\\
\les&  \epsilon^{2-\kappa}\sup_{L\leq \frac{\epsilon d}2}L^{-\kappa}\inf_l\|U_{T}(x,.)-l\|_{B(x,L)}\\
&+ \sum_{\beta\in A}M^{(1)}_{D_d,\frac{d}2}\epsilon^{-4+2\beta-\kappa}+\epsilon^{2-2\kappa}\frac{ d^{-\kappa}}{2^{-\kappa}}\|U(x,.)\|_{B(x,\frac{d}2(1+\epsilon^2))},
\end{align*}
where the ratios between $L$ and $T$, and $R$ and $T$ are fixed only within the supremum operators.
Applying Equation \eqref{Schauder Step 2} gives
\begin{align}\label{app_st_21}
\sup_{R\leq \frac{\epsilon d}2}R^{-\kappa}&\inf_l\|U(x,.)-l\|_{B(x,R)}\nonumber\\
\leq\sup_{R\leq \frac{\epsilon d}2}&R^{-\kappa}\inf_l\|U_{T}(x,.)-l\|_{B(x,R)}+M^{(2)}_{\{x\},\frac{\epsilon d}2,\frac{\epsilon^2 d}2} \sum_{\beta\in A}\epsilon^{\kappa-\beta}
+ \epsilon^\kappa[U]_{\kappa, B(x,\frac{\epsilon d}2),\frac{\epsilon^2 d}2}\nonumber\\
\les \sum_{\beta\in A}&\Big(M^{(1)}_{\{x\},\frac{d}2}\epsilon^{-4+2\beta-\kappa}+M^{(2)}_{D_d,\frac{\epsilon d}2,\frac{\epsilon^2 d}2} \epsilon^{\kappa-\beta}\Big)+\epsilon^{2-2\kappa}\frac{ d^{-\kappa}}{2^{-\kappa}}\|U(x,.)\|_{B(x,\frac{d}2(1+\epsilon^2))}\nonumber\\
+ \epsilon^\kappa&[U]_{\kappa, B(x,\frac{\epsilon d}2),\frac{\epsilon^2 d}2}+ \epsilon^{2-\kappa}\sup_{L\leq \frac{\epsilon d}2}L^{-\kappa}\inf_l\|U_{T}(x,.)-l\|_{B(x,L)}\\
\les \sum_{\beta\in A}&\Big(M^{(1)}_{\{x\},\frac{d}2}\epsilon^{-4+2\beta-\kappa}+M^{(2)}_{\{x\},\frac{\epsilon d}2,\frac{\epsilon^2 d}2} \epsilon^{\kappa-\beta}+M^{(2)}_{\{x\},\frac{ d}2,\frac{\epsilon^2 d}2} \epsilon^{2(\kappa-\beta)}\Big)\nonumber\\
&+\epsilon^{2-2\kappa}\frac{ d^{-\kappa}}{2^{-\kappa}}\|U(x,.)\|_{B(x,\frac{d}2(1+\epsilon^2))}+ (\epsilon^\kappa+\epsilon^{2+\kappa})[U]_{\kappa, B(x,\frac{\epsilon d}2),\frac{\epsilon^2 d}2}\nonumber\\
&+ \epsilon^{2-\kappa}\sup_{L\leq \frac{\epsilon d}2}L^{-\kappa}\inf_l\|U(x,.)-l\|_{B(x,L)}.\nonumber
\end{align}
The last term on the \rhs can now be absorbed into the \lhs for $\epsilon$ sufficiently small, giving the bound \ref{app_st_2}

\paragraph*{\textsc{Step} 4.}
 We prove that 
\begin{align}\label{Schauder Step 3}
\sup_{d\leq d_0}d^\kappa&[U]_{\kappa,D_d}\nonumber\\
\les &\sum_{\beta\in A}\Big(M^{(1)}\epsilon^{-4+2\beta-\kappa}+M^{(2)} \epsilon^{\kappa-\beta}\Big)+(\epsilon^{-\kappa}+\epsilon^{2-2\kappa})\sup_{d\leq d_0}\|U\|_{D_d,d}.
\end{align}
We first argue that we can change the order of the supremum and the infimum in $\sup_{R\leq \frac{\epsilon d}2}R^{-\kappa}\inf_l\|U(x,.)-l\|_{B(x,R)}$.
 Since $U(x,x)=0$ it is clear that one can restrict to $l(x)=0$ hence $l(y)=C(x,R).X(y-x)$.
  We argue that $C$ may be chosen independently of $R$. Let $C_R$ be the (near) optimal constant for the radius $R$. Then 
 \[
 R^{-(\kappa-1)}|C_{\frac{R}2}-C_R|\lesssim \sup_{R\leq \frac{\epsilon d}2}R^{-\kappa}\inf_l\|U(x,.)-l\|_{B(x,R)}.
 \]
  Since $\kappa>1$, this can be extended by summation to all $R\leq \frac{\epsilon d}2$, thus there exists a near optimal constant $C$ independent of $\rho$. We then have
\begin{align*}
\inf_{\nu(x)}\sup_{x\neq y\in B(x,\frac{\epsilon d}2)}d(x,y)^{-\kappa}|U(x,y)-\nu(x).X(y-x)|\\
\leq\inf_l\sup_{R\leq \frac{\epsilon d}2}R^{-\kappa}\|U(x,.)-l\|_{B(x,R)}\\
\les\sup_{R\leq \frac{\epsilon d}2}R^{-\kappa}\inf_l\|U(x,.)-l\|_{B(x,R)}.
\end{align*}
Therefore, if we take the supremum over $x\in D_d$ in Equation \eqref{app_st_2} then multiply it by $d^\kappa$ and take the supremum over $d$, we get
\begin{align*}
\sup_{d\leq d_0}&d^\kappa\sup_{x\in D_d}\inf_{\nu(x)}\sup_{x\neq y\in B(x,\frac{\epsilon d}2)}d(x,y)^{-\kappa}|U(x,y)-\nu(x).X(y-x)|\\
\les &\sup_{d\leq d_0}d^\kappa\sum_{\beta\in A}\Big(M^{(1)}_{D_d,\frac{d}2}\epsilon^{-4+2\beta-\kappa}+M^{(2)}_{D_d,\frac{ d}2,\frac{\epsilon^2 d}2} \epsilon^{\kappa-\beta}\Big)+\epsilon^{2-2\kappa}\sup_{d\leq d_0}\|U\|_{D_d,d}\nonumber\\
+&\epsilon^\kappa\sup_{d\leq d_0}d^\kappa\sup_{x\in D_d}\sup_{y\in B(x,\frac{\epsilon d}2)}\inf_{\nu(y)}\sup_{y\neq z\in B(y,\frac{\epsilon^2 d}2)}d(y,z)^{-\kappa}|U(y,z)-\nu(y).X(z-y)|.
\end{align*}
The last term can be absorbed into the left-hand side for $\epsilon$ small enough since for $y\in B(x,\frac{\epsilon d}2)$ we have $d(y,\delta D)\geqslant d(1-\frac\epsilon2)$ and consequently for $z\in B(y,\frac{\epsilon^2 d}2)$, we have $d(y,z)\leq \frac{\epsilon^2d(y,\delta D)}{2(1-\frac\epsilon2)}\leq \frac{\epsilon d(y,\delta D)}2$, which gives
\begin{align*}
\sup_{d\leq d_0}&d^\kappa\sup_{x\in D_d}\inf_{\nu(x)}\sup_{x\neq y\in B(x,\frac{\epsilon d}2)}d(x,y)^{-\kappa}|U(x,y)-\nu(x).X(y-x)|\\
\les &\sup_{d\leq d_0} d^\kappa\sum_{\beta\in A}\Big(M^{(1)}_{D_d,\frac{d}2}\epsilon^{-4+2\beta-\kappa}+M^{(2)}_{D_d,\frac{\epsilon d}2,\frac{\epsilon^2 d}2} \epsilon^{\kappa-\beta}\Big)+\epsilon^{2-2\kappa}\sup_{d\leq d_0}\|U\|_{D_d,d}.
\end{align*}
We concludes the proof of \eqref{Schauder Step 3} by extending to all $y\in D_d$ with the following argument
\begin{align*}
\sup_{x\in D_d}&\inf_{\nu(x)}\sup_{x\neq y\in D_d}d(x,y)^{-\kappa}|U(x,y)-\nu(x).X(y-x)|\\
\leq &\sup_{x\in D_d}\inf_{\nu(x)}\sup_{x\neq y\in B(x,\frac{\epsilon d}2)}d(x,y)^{-\kappa}|U(x,y)-\nu(x).X(y-x)|\\
&+\sup_{x\in D_d}\inf_{\nu(x)}\sup_{y\in D_d\setminus B(x,\frac{\epsilon d}2)}d(x,y)^{-\kappa}|U(x,y)-\nu(x).X(y-x)|\\
\leq& \sup_{x\in D_d}\inf_{\nu(x)}\sup_{x\neq y\in B(x,\frac{\epsilon d}2)}d(x,y)^{-\kappa}|U(x,y)-\nu(x).X(y-x)|+\Big(\frac{\epsilon d}2\Big)^{-\kappa}\|U\|_{D_d,d}.
\end{align*}

\subsubsection{Proof of Corollary \ref{corschauder}} 
From the definition of $[U]_{\kappa,D}$ in \eqref{e:def-hol2var} used with variables $x,y\in D_d$ and with triangle inequalities, we get
\[
|\nu(x).X(y-x)|\leq [U]_{\kappa,D_d}d(x,y)^\kappa+\|U\|_{D_d,d(x,y)}.
\]
Applying the interior cone condition for $r\in [0,r_d]$ gives the existence of some $y$ with $d(x,y)=r$ such that
\[
\lambda|\nu(x)|d(x,y)\leq [U]_{\kappa,D_d}d(x,y)^\kappa+\|U\|_{D_d,r},
\]
which proves \eqref{corschauder1}.

Using again the definition of $[U]_{\kappa,D}$ with variables $x,y$ and $y,z\in D_d$, and with triangle inequalities, we get
\begin{align*}
|U(x,y)-U(x,z&)-U(z,y)-(\nu(x)-\nu(z)).X(y-z)|\\
&\leq [U]_{\kappa,D_d}(d(x,y)^\kappa+d(y,z)^\kappa+d(x,z)^\kappa).
\end{align*} 
We combine this with the three-point continuity condition \eqref{lemschauder2}, and we assume that $r>d(x,z)=d(y,z)\geqslant \frac{d(x,y)}2$ to get
\[
|(\nu(x)-\nu(z)).X(y-z)|\les d(x,z)^\kappa([U]_{\kappa,D}+M^{(2)}_{D_d,\frac{d}4,\frac{d}4}).
\]
Choosing finally $y$ such that $|(\nu(x)-\nu(z)).X(y-z)|\geqslant \lambda|\nu(x)-\nu(z)||d(y,z)|$ gives \eqref{corschauder2} for $d(x,y)\leq r$. For $d(x,y)\geqslant r$, we have
\[
d(x,y)^{-\kappa+1}|\nu(x)-\nu(y)|\leq 2  r^{-\kappa+1}\|\nu\|_{D_d}.
\]
Applying \eqref{corschauder1} gives for $d(x,y)\geqslant r$,
\[
d(x,y)^{-\kappa+1}|\nu(x)-\nu(y)|\les [U]_{\kappa,D_d}+ r^{-\kappa}\|U\|_{D_d,r}.
\]
\subsection{Proof of the reconstruction}\label{proof reconstruction}
In this section we will use the following notations for $f$ a function of one variable and  $F$ a function of two variables:
\begin{align}\label{notations}
%[f,(\cdot)_T](x)&=f(x)-f_T(x)=\int\Psi_T(x-y)(f(x)-f(y))dy\nonumber\\
[F,(\cdot)_T](x)&=\int\Psi_T(x-y)F(x,y)dy\\
[Ff,(\cdot)_T](x)&=\int\Psi_T(x-y)F(x,y)f(y)dy\nonumber.
\end{align}
\subsubsection{Proof of Theorem \ref{Reconstruction}}
This is the only place where our particular choice of convolution kernel is crucial. It allows to use the following factorisation:
\begin{align*}
\Big|[F,(\cdot)_{T2^{-n}}]&(x_1)-\Big([F,(\cdot)_{T2^{-n-1}}]\Big)_{T2^{-n},1}(x_1)\Big|\\
=&\Big|\int\int\Psi_{T2^{-n-1}}(x_2-y)\Phi_{T2^{-n-1}}(x_1-x_2)(F(x_1,y)-F(x_2,y))dydx_2\Big|\\
\leq&\sum_{\beta\in A}C_\beta\int\Phi_{T2^{-n-1}}(x_1-x_2)d(x_1,x_2)^{\gamma_\beta-\beta}(T2^{-n-1})^\beta dx_2\\
\leq&\sum_{\beta\in A}C_\beta(T2^{-n-1})^{\gamma_\beta}.
\end{align*}
This proves the convergence of $[F,(\cdot)_{T2^{-n}}]$ to $f:y\mapsto F(y,y)$ and justifies the bound following telescopic sum, obtained once more thanks to the semi-group property of our kernel:
\begin{align*}
\Big|[F,(\cdot)_T]-[F,(\cdot)_{T2^{-N}}]_{T,N-1}\Big|
=&\Big|\sum_{n=0}^N\Big([F,(\cdot)_{T2^{-n}}]-[F,(\cdot)_{T2^{-n-1}}]_{T2^{-n},1}\Big)_{T,n}\Big|\\
\leq&\sum_{n=0}^N\sum_{\beta\in A}C_\beta(T2^{-n-1})^{\gamma_\beta}\les \sum_{\beta\in A}C_\beta T^{\gamma_\beta},
\end{align*}
where the constant in ''$\les$'' depends only on $\gamma$ (in particular not on $N$), thus proving the theorem.

\subsubsection{Proof of Lemma \ref{lem_v2ou}}
To obtain a bound on $(v^2\ou)_T(y)$ we implement the following expansion.
\begin{align*}
(v^2\ou)_T(x)=
(v^2\ou)_T(x)&-v(x)^2\ou_T(x)-2v(x)((\outz(x)-\outz)\ou)_T(x)\\
&+v(x)^2\ou_T(x)+2v(x)((\outz(x)-\outz)\ou)_T(x).
\end{align*}
From the bound \eqref{def_outu} we have 
\begin{align}
|v(x)^2\ou_T(x)|\leq \|v\|_{B(x,T)}^2[\ou]_{-\frac12-\epsilon}T^{-\frac12-\epsilon}\label{chose 0}\\
|v(x)((\outz(x)-\outz)\ou)_T(x)|\leq \|v\|_{B(x,T)}[\outu]_{-4\epsilon}T^{-4\epsilon}. \label{chose 2}
\end{align}
To bound the remaining part, we will apply Theorem~\ref{Reconstruction} and to that end we set 

\begin{equation}\label{defFv2ou}
F(x_1,y)=v(x_1)^2\ou(y)+2v(x_1)(\outz(x_1)-\outz(y))\ou(y).
\end{equation}
Then 
\begin{align*}
F(x_1,y)-F(x_2,y)=&(v(x_1)+v(x_2))(v(x_1)-v(x_2)+\outz(x_1)-\outz(x_2))\ou(y)\\
&+(v(x_1)-v(x_2))(\outz(x_1)-\outz(x_2))\ou(y)\\
&+2(v(x_1)-v(x_2))(\outz(x_1)-\outz(y))\ou(y).
\end{align*}
By definition of $\outu$ \eqref{def_outu} this gives, for $x_1,x_2\in B(x,T-t)$
\begin{align*}
|\int\Psi_t(x_1-y)&(F(x_1,y)-F(x_2,y))dy|\\
\leq& 2\|v\|_{D_d}[v+\outz]_{1-2\epsilon,B(x,T)}d(x_1,x_2)^{1-2\epsilon}[\ou]_{-\frac12-\epsilon}t^{-\frac12-\epsilon}\\
&+[v]_{\frac12-3\epsilon,B(x,T)}[\outz]_{\frac12-3\epsilon}d(x_1,x_2)^{1-6\epsilon}[\ou]_{-\frac12-\epsilon}t^{-\frac12-\epsilon}\\
&+[v]_{\frac12-3\epsilon,B(x,T)}d(x_1,x_2)^{\frac12-3\epsilon}[\outu]_{-4\epsilon}t^{-4\epsilon}.
\end{align*}
Hence by Theorem~\ref{Reconstruction}, we have the bound
\begin{align}\label{app_rec_v2ou}
|(v^2\ou)_T(x)-v(x)^2\ou_T(x)+2v(x)((\outz(x)-\outz)\ou)_T(x)	|\nonumber\\
\les T^{\frac12-3\epsilon}2\|v\|_{D_d}[v+\outz]_{1-2\epsilon,B(x,T)}[\ou]_{-\frac12-\epsilon}\nonumber\\
+T^{\frac12-7\epsilon}[v]_{\frac12-3\epsilon,B(x,T)}([\outz]_{\frac12-3\epsilon}[\ou]_{-\frac12-\epsilon}+[\outu]_{-4\epsilon}).
\end{align}
Together with bounds \eqref{chose 0} to \eqref{chose 2}, we get lemma \ref{lem_v2ou}.

\subsubsection{Proof of Lemma \ref{lem_voud}}

The last quantity we need to bound is:
\begin{align*}
((v-v(x))\oud)_T(x)+&3C_2(v_T(x)+\ou_T(x))=\\
&[\bar{U}\oud,(\cdot)_T](x)-3C_2(v-v_T)(x)-((\outz-\outz(x))\oud-3C_2\ou)_T(x)\\
&-3v(x)((\oudz-\oudz(x))\oud-3C_2)_T(x)+\nu(x).(X(\cdot-x)\oud)_T(x),
\end{align*}
where $\nu$ is optimal in the definition of $[U]_{1+3\epsilon,D}$ and $\bar{U}(x,y)=U(x,y)-\nu(x).X(y-x)$.
From the bounds \eqref{def_outd}, \eqref{def_oudd} and \eqref{def_oudX} we have 
\begin{align}
|((\outz-\outz(x))\oud-3C_2\ou)_T(x)|\leq& [\outd]_{-\frac12-5\epsilon}T^{-\frac12-5\epsilon}\label{chose 1},\\
|v(x)((\oudz-\oudz(x))\oud-3C_2)_T(y)\leq& \|v\|_{B(x,T)}[\oudd]_{-4\epsilon}T^{-4\epsilon},\\
|\nu(x).(X(\cdot-y)\oud)_T(x)|\leq& \|\nu\|_{B(x,T)}[\oudX]_{-2\epsilon}T^{-2\epsilon}.\label{chose 3}
\end{align}
To bound the remaining part, we will apply Theorem~\ref{Reconstruction} and to that end we set 
\begin{align}\label{defFvoud}
F(x_1,y)=(v(x_1)+\outz(x_1)-\outz(y)+3v(x_1)(\oudz(x_1)-\oudz(y))-\nu(x_1).X(x_1-y))\oud(y)\nonumber\\-3C_2(v(x_1)-v(y)).
\end{align}
Then for $x_1,x_2\in B(x,T-t)$,
\begin{align*}
F(x_1,y)-F(x_2,y)=&(3(v(x_1)-v(x_2))((\oudz(y)-\oudz(x_2))\oud(y)-C_2)\\
&+U(x_1,x_2)\oud(y)-(\nu(x_1)-\nu(x_2)).X(y-x_2)\oud(y).
\end{align*}
By definition of $\oudd$ \eqref{def_oudd} and $\oudX$ \eqref{def_oudX}, this gives
\begin{align*}
\int\Psi_t(x_1-y)&(F(x_1,y)-F(x_2,y))dy\\
\leq&3[v]_{\frac12-3\epsilon,B(x,T)}d(x_1,x_2)^{\frac12-3\epsilon}[\oudd]_{-4\epsilon}t^{-4\epsilon} \\
&+[U]_{\frac32-5\epsilon,B(x,T)}d(x_1,x_2)^{\frac32-5\epsilon}[\oud]_{-1-2\epsilon}t^{-1-2\epsilon}\\
&+[\nu]_{\frac12-5\epsilon,B(x,T)}d(x_1,x_2)^{\frac12-5\epsilon}[\oudX]_{-2\epsilon}t^{-2\epsilon}.
\end{align*}
Hence by Theorem~\ref{Reconstruction}, we have the bound
  $v\oud$ such that 
\begin{align}\label{app_rec_voud}
\|[\bar{U}\oud,(\cdot)_T]&(x)-3C_2(v-v_T)(x)\|\nonumber\\ 
\les T^{\frac12-7\epsilon}&\Big(3[v]_{\frac12-3\epsilon,B(x,T)}[\oudd]_{-4\epsilon}
+[U]_{\frac32-5\epsilon,B(x,T)} [\oud]_{-1-2\epsilon}+[\nu]_{\frac12-5\epsilon,B(x,T)}[\oudX]_{-2\epsilon}\Big).
\end{align}
Together with bounds \eqref{chose 1} to \eqref{chose 3}, we get the lemma \ref{lem_voud}.

\subsection{Proof of Lemma \ref{main thm poly 1}}\label{proof max}
This proof is a version of the proof of Theorem 4.4 in our companion paper, \cite{2018arXiv180810401M}, specialised to cubic
non-linearity.  This specialisation makes the proof significantly simpler.
We only prove the bound for the positive part of $u$. The bound for the negative part follows by symmetry.
Let $\eta$ be a continuous function defined on $\R_+\times[-1,1]^3$, $C^2$ and strictly positive on the interior and such that $\eta=0$ on the boundary. 
Either $u\eta$ attains its maximum on $[0,1]\times [-1,1]^3$ at some point $z_0\in(0,1]\times (-1,1)^3$, or it is non-positive, in which case $u\leqslant0$ in $[0,1]\times [-1,1]^3$. Assuming this is not the case, we get that at the maximum point, $0=\grad (u\eta)(z_0)$, i.e. 
\begin{equation}\label{eq:max grad u eta}
\grad u=-\frac{\grad\eta}{\eta}u.
\end{equation}
If $z_0\in \{1\}\times (-1,1)^3$, then $\partial_tu\eta(z_0)\geqslant 0$. Else, $\partial_tu\eta(z_0)= 0$. Additionally, $\Delta u\eta(z_0)\leqslant 0$ and therefore at the maximum we have
\begin{align*}
0\leqslant& \heat(u\eta)=\eta\heat u+u\heat\eta-2\grad u.\grad \eta\\
\overset{\eqref{eq:max rd 1};\eqref{eq:max grad u eta}}{=} &-\eta (u^3-g(u,z))+u\Big(\heat\eta+2\frac{|\grad\eta|^2}{\eta}\Big).
\end{align*}
Assume $\eta$ satisfies the following inequality:
\begin{equation}\label{eq:max eta}
 \frac{\heat\eta}\eta+2\frac{|\grad\eta|^2}{\eta^2}\leqslant \frac1{2\eta^2}.
\end{equation}
Then we get
\begin{equation}\label{eq:max u eta}
 u^2\leqslant \frac{\eta^{-2}}2+\frac{\|g\|}{u}\leqslant 2\max\Big{\{}\frac{\eta^{-2}}2,\frac{\|g\|}{u}\Big{\}}.
\end{equation}
If the maximum on the \rhs is realised by the first term, then at $z_0$, $u\eta\leqslant 1$. If the maximum is realised by the second term, then it has to be bigger than the first one:
\[
\frac{\eta^{-2}}2\leqslant\frac{\|g\|}{u}\Rightarrow u\eta\leqslant2\eta^3\|g\|.
\]
We then have that at $z_0$, $u\eta\leqslant 2$ under the condition
$\eta\leqslant\|g\|^{-\frac13}$. In both cases, we obtain that $u\leqslant \frac2\eta$ on all of $(0,1]\times (-1,1)^3$. With a choice of $\eta$  also satisfying the inequality \eqref{eq:max eta}, we obtain good bounds on the function $u$. We choose the following for $z=(t,x)\in(0,\infty)\times(-1,1)^3$, for some value $\lambda$ to be defined:
\begin{equation}\label{eq:max def eta}
\eta(x,t)=\frac{\lambda}{\lambda\|g\|^\frac13+\frac{1}{\sqrt{t}}+\sum_{i=1}^d\frac{1}{1+x_i}+\frac{1}{1-x_i}},
\end{equation}
and we continuously extend with the value $0$ on the boundary of the domain.
This choice of $\eta$ guarantees a bound on $u$ that is related to the distance from the boundary of $[0,1]\times[-1,1]^d$, independently of the boundary conditions. Indeed,
\begin{align}\label{eq:max comp eta}
( 2d+1)\frac1{\lambda\min_i\{ \sqrt{t},1+x_i,1-x_i\}}\geqslant\frac1\eta-\|g\|^\frac13\geqslant\frac1{\lambda\min_i\{ \sqrt{t},1+x_i,1-x_i\}}.
\end{align}
It also satisfies $0\leqslant\eta\leqslant\|g\|^{-\frac13}$.
 We compute the derivatives to check \eqref{eq:max eta}.
\[
\partial_t\eta=\frac1{2\lambda t^{\frac32}}\eta^2,\qquad \partial_i\eta=\frac{1}{\lambda}\Big(\frac{1}{(1+x_i)^2}-\frac{1}{(1-x_i)^2}\Big)\eta^2,
\]
and
\[
\partial_i^2\eta=-\frac{2x_i}{\lambda}\Big(\frac{1}{(1+x_i)^3}+\frac{1}{(1-x_i)^3}\Big)\eta^2
+\frac{2}{\lambda^2}\Big(\frac{1}{(1+x_i)^2}-\frac{1}{(1-x_i)^2}\Big)^2\eta^3.
\]
Some of these terms cancel and we have the bound, using \eqref{eq:max comp eta}
\begin{align*}
2\eta\heat\eta+4|\grad\eta|^2=&\frac1{\lambda t^{\frac32}}\eta^3+4\sum_{i=1}^3\frac{x_i}{\lambda}\Big(\frac{1}{(1+x_i)^3}+\frac{1}{(1-x_i)^3}\Big)\eta^3
\leq25\lambda^2.
\end{align*}
Therefore, by taking $\lambda\leq \frac15$, we have proved Lemma \ref{main thm poly 1}.

\bibliographystyle{abbrv}
\bibliography{phip}

\end{document}